\documentclass[11pt]{article}
\usepackage{graphicx}
\usepackage{amsfonts}
\usepackage{amsmath}
\usepackage{amssymb}
\usepackage{fancyhdr}
\usepackage{titlesec}
\usepackage{indentfirst}
\usepackage{booktabs}
\usepackage{verbatim}
\usepackage{color}
\usepackage{amsthm}
\usepackage[normalem]{ulem}
\usepackage{hyperref}

\usepackage{bm}
\usepackage{pdflscape}

\usepackage{todonotes}

\usepackage{tikz}
\usetikzlibrary{arrows}
\usetikzlibrary{decorations.markings}
\usetikzlibrary{pgfplots.groupplots}
\usepackage{wrapfig}
\usepackage{pgfplots}
\pgfplotsset{compat=1.18}
\usepackage{cleveref}
\usepackage{algorithm}
\usepackage{algorithmic}
\usepgfplotslibrary{fillbetween}

\usepackage{cite} %

\usepackage[page,header]{appendix}
\usepackage{titletoc}

\newcommand{\rd}{\,\mathrm{d}}

\usepackage{multirow}

\newcommand{\mathR}{\mathbb{R}}

\newcommand{\bG}{\mathbf{G}}

\newcommand{\bc}{\mathbf{c}}
\newcommand{\bU}{\mathbf{U}}
\newcommand{\bR}{\mathbf{R}}
\newcommand{\bI}{\mathbf{I}}

\newcommand{\bmu}{\bm{\mu}}

\newcommand{\bbeta}{\bm{\beta}}
\newcommand{\bx}{\mathbf{x}}
\newcommand{\bv}{\mathbf{v}}
\newcommand{\bu}{{u}}
\newcommand{\bn}{\mathbf{n}}
\newcommand{\bbf}{{f}}
\newcommand{\bbg}{{g}}

\newcommand{\Bbar}{\bar{B}}

\newcommand{\rbF}{\boldsymbol{\mathcal{F}}}
\newcommand{\mL}{\mathcal{L}}

\newcommand{\RR}{\mathbb{R}}
\newcommand{\mathV}{\mathbb{V}}
\newcommand{\Nt}{\mathcal{N}}

\newcommand{\fFOM}{\bbf}

\newcommand{\ffROM}{\bbf_{\text{rb}}}
\newcommand{\fflat}{\mathbf{f}}
\newcommand{\gflat}{\mathbf{g}}

\newcommand{\Res}{\text{Res}}
\newcommand{\MinRes}{\text{MinRes}}

\newcommand{\mM}{\mathcal{M}}
\newcommand{\Mtrain}{\mM_{\text{train}}}
\newcommand{\Mtest}{\mM_{\text{test}}}
\newcommand{\Mrange}{\mM_{\text{range}}}
\newcommand{\Meim}{\mM_{\text{EIM}}}
\newcommand{\mueim}{\tilde{\bmu}}

\newcommand{\Tdata}{T^{\text{data}}}
\newcommand{\mass}{m_0}

\newtheorem{theorem}{Theorem}[section]

\newtheorem{remark}[theorem]{Remark}

\definecolor{matlabblue}{RGB}{0,114,189}
\definecolor{matlaborange}{RGB}{217,83,25}
\definecolor{matlabyellow}{RGB}{237,177,32}
\definecolor{matlabpurple}{RGB}{126,47,142}
\definecolor{matlabgreen}{RGB}{119,172,48}
\definecolor{matlabcyan}{RGB}{77,190,238}

\pgfplotscreateplotcyclelist{siam}{
    {matlabblue, thick, solid, mark=square*, mark size=1.8pt, mark options={solid}},
    {matlaborange, thick, densely dashdotted, mark=*, mark size=1.8pt, mark options={solid}},
    {matlabyellow, thick, densely dashed, mark=triangle*, mark size=1.8pt, mark options={solid}},
    {matlabpurple, thick, densely dotted, mark=diamond*, mark size=1.8pt, mark options={solid}},
    {matlabgreen, thick, dashdotted, mark=pentagon*, mark size=1.8pt, mark options={solid}},
    {matlabcyan, thick, dashdotdotted, mark=star, mark size=1.8pt, mark options={solid}},
}

\begin{document}

	\title{Reduced order model for parametric Boltzmann equation and its application to inverse problems}
	\author{
    Shanyin Tong\thanks{tong3@sas.upenn.edu. Department of Mathematics, University of Pennsylvania, Philadelphia, PA 19104}
\and
Jingwei Hu\thanks{hujw@uw.edu. Department of Applied Mathematics, University of Washington, Seattle, WA 98195}
\and Fengyan Li\thanks{lif@rpi.edu. Department of Mathematical Sciences, Rensselaer Polytechnic Institute, Troy, NY 12180}
\and
Zhaiming Shen\thanks{zshen49@gatech.edu.
School of Mathematics,
Georgia Institute of Technology, Atlanta, GA 30332}
\and Yunan Yang\thanks{yunan.yang@cornell.edu. Department of Mathematics, Cornell University, Ithaca, NY 14850}
\and
Jiajia Yu\thanks{yujiajia.math@gmail.com. Department of Mathematics, Duke University, Durham, NC 27710}
}
	\date{}
	
	\maketitle

\begin{abstract}
    The Boltzmann equation plays an important role in modeling mesoscopic behavior in a wide range of scientific and engineering applications. However, its numerical solution is computationally expensive due to the high dimensionality of the model and the nonlinear nonlocal collision operator, especially for steady-state problems that require iterative solvers. This cost becomes prohibitive for inverse problems, where the induced optimization problem requires repeated forward solves. 
    
    In this work, we propose a reduced-order model (ROM) for the parametric Boltzmann equation to address this computational challenge. The ROM constructs a low-dimensional approximation space for the parameter-induced  solution manifold through a residual-based greedy strategy, and
    the reduced solution is then obtained via residual minimization over the reduced space,
    subject to mass conservation. The overall efficiency of the ROM is achieved by exploiting the quadratic structure of the collision operator and a precomputed separable approximation of the collision kernel. 
    The resulting ROM is further applied to a thermally-driven inverse problem  for reconstructing collision parameters from the observed macroscopic temperature
data.
This is accomplished either by directly replacing the PDE constraint with the ROM, leading to a bilevel optimization formulation, or by reformulating the task
    as a single-level optimization problem through the Karush--Kuhn--Tucker (KKT) conditions.
Numerical experiments in both collision-dominated and transport-dominated cases are performed to demonstrate the efficiency and accuracy of the proposed ROM and its effectiveness in inverse problems. In particular, the resulting inverse problem is computationally much more tractable, achieving speedups of several orders of magnitude 
over that based on the full-order model while maintaining comparable accuracy.
\end{abstract}

\section{Introduction}

The Boltzmann equation \cite{Cercignani} describes the non-equilibrium dynamics of complex particle systems through a probability density distribution, and has found applications in many scientific and engineering disciplines, including rarefied gas dynamics, plasma physics, and even the biological and social sciences. Despite providing a more accurate description than fluid models such as the Navier--Stokes equations, numerical solutions of the Boltzmann equation require resolving the high-dimensional phase space (e.g., up to six dimensions) and approximating the nonlinear collision integral operator. In practice, many problems of interest also require solving the steady-state Boltzmann equation, for which iterative methods often require anywhere from hundreds to hundreds of thousands of iterations to converge. Taken together, these factors make solving the Boltzmann equation extremely computationally demanding.

The situation is further exacerbated when parameter studies or inverse problems are required. Indeed, the collision kernel depends on the underlying physical system and is often described by a parametric form. Studying the dependence of the solution on these parameters or inferring them from observational data makes the task even more formidable, given the already high computational cost of full-order forward simulations.

Motivated by these challenges, in this work, we develop a model order reduction (MOR) strategy that provides a reduced-order model 
(ROM) for computing steady-state solutions of the parametric Boltzmann equation at significantly lower computational complexity, when compared with standard grid-based numerical methods, for a large collection of parameter values. This reduced-order surrogate solver is further employed in an inverse problem framework  to infer the parameters in the collisional kernel from observational data  at  a practically manageable cost. 

With the ROM, the complexity reduction of the forward problem is achieved by integrating several technical  ingredients: a greedy strategy as in reduced basis methods \cite{hesthaven2016certified, haasdonk2017reduced} to build a low-dimensional approximation/representation space for the parameter-induced solution manifold, from which a reduced solution is sought via residual minimization; full exploitation  of the intrinsic quadratic structure of the collision operator combined with a precomputed separable approximation of the parametric collision kernel.  The inverse problem is solved through its reformulation into  either a bilevel or a single-level optimization problem motivated by algorithms proposed in~\cite{guerra2025learning}. Our method efficiently harvests from the aforementioned low-dimensional representation of the solution and applies to the parametric quadratic collision operator.

Recently, there has been growing interest in developing MOR techniques for kinetic models and inverse problems. We first review related work on MOR for kinetic models. A growing body of work has emerged, particularly for linear steady-state kinetic models such as radiative transfer and neutron transport, using strategies including proper orthogonal decomposition (POD) \cite{behne2022minimally, behne2023parametric}, proper generalized decomposition (PGD) \cite{dominesey2022reduced}, and the reduced basis method (RBM)  \cite{peng2022reduced, matsuda2025reduced}. By contrast, relatively little progress has been made on the nonlinear Boltzmann equation, particularly for the computation of steady-state solutions.

We next turn to some work on MOR for inverse problems. Inverse problems often rely on forward solvers, such as numerical methods for the underlying PDE models, whose repeated evaluations can be computationally expensive if not prohibitive. A natural strategy to reduce the overall cost of these multi-query tasks is to replace the expensive forward solver with an efficient surrogate model. Examples of ROMs applied to inverse problems include
reconstructing Arrhenius parameters in combustion models \cite{galbally2010non}, reconstructing conductivity fields in groundwater flow problems \cite{lieberman2010parameter},  full waveform inversion \cite{borcea2024data}, and parameter estimation in linear or semilinear parabolic PDEs \cite{schmidt2013derivative, kartmann2025adaptive}.

Before presenting the organization and main text of the paper, we highlight the main contributions of this work:
(1) To the best of our knowledge, this is the first attempt to develop a reduced-order model for the parametric steady-state nonlinear Boltzmann equation. We directly exploit the intrinsic quadratic structure of the collision operator, together with an efficient separable approximation of the collision kernel using EIM. These two ingredients form
the basis for the computational acceleration of the proposed ROM.
(2) Beyond accelerating the forward problem, we demonstrate the effectiveness and efficiency of the proposed ROM for a challenging inverse problem that reconstructs microscopic collision parameters from macroscopic observations. We further propose a single-level reformulation of the inverse problem based on the reduced basis space, which reduces the computational cost by several thousand times while maintaining reconstruction accuracy comparable to that of the full-order model.
(3) Numerical experiments demonstrate the applicability and robustness of the proposed framework in both collision-dominated and transport-dominated regimes.

The rest of the paper is organized as follows. In \Cref{sec:Boltz}, we introduce the mathematical formulation of the Boltzmann equation and its inverse problem. In \Cref{sec:FOM}, we discuss the numerical schemes for solving the Boltzmann equation, which we later refer to as the full-order model (FOM). In \Cref{sec:ROM}, we present the proposed reduced-order model for the Boltzmann equation. We first explain the efficient evaluation of the collision operator using its quadratic structure, EIM approximation and precomputation of the collision kernel, and then describe the resulting greedy offline construction and online evaluation schemes. In \Cref{sec:IVP}, we discuss the inverse problem and how the ROM can be applied through either a bilevel or a single-level optimization formulation. Finally, in \Cref{sec:num}, we present numerical experiments to demonstrate the efficiency and accuracy of the ROM itself and its effectiveness in  an inverse problem.
The paper is concluded in \Cref{sec:con}.

\section{The Boltzmann equation and its inverse problem}\label{sec:Boltz}	

We consider the time-independent Boltzmann equation:
\begin{equation} \label{eq:full_Boltz}
\bv\cdot\nabla_{\bx}f(\bx,\bv)=Q(f,f)(\bx,\bv), \quad \bx\in \Omega\subset \mathbb{R}^{d}, \quad \bv\in \mathbb{R}^{d}, \quad d=2 \text{ or }3,
\end{equation}
where $f=f(\bx,\bv)$ is the phase-space distribution function in position $\bx$ and velocity $\bv$, and $Q(f,f)$ is the Boltzmann collision operator given by
\begin{equation} \label{eq:full_collision}
Q(f,f)(\bx,\bv)=\int_{\mathbb{R}^{d}}\int_{\mathbb{S}^{d-1}} B(\bv-\bv_*,\sigma)[f(\bx,\bv_*')f(\bx,\bv')-f(\bx,\bv_*)f(\bx,\bv)]\rd{\sigma}\rd{\bv_*}.
\end{equation}
Here, $(\bv,\bv_*)$ and $(\bv',\bv_*')$ are the velocity pairs before and after a collision, related by conservation of momentum and energy:
\begin{equation}
\bv'=\frac{\bv+\bv_*}{2}+\frac{|\bv-\bv_*|}{2}\sigma, \quad \bv_*'=\frac{\bv+\bv_*}{2}-\frac{|\bv-\bv_*|}{2}\sigma,
\end{equation}
with $\sigma$ varying on the unit sphere $\mathbb{S}^{d-1}$. The function $B$ is the collision kernel that determines the collision frequency; it depends only on $|\bv-\bv_*|$ and on the cosine of the deviation/scattering angle:
\begin{equation}
\label{eq:kernel}
B(\bv-\bv_*,\sigma)=B(|\bv-\bv_*|,\cos\theta), \quad \cos\theta=\sigma \cdot \frac{\bv-\bv_*}{|\bv-\bv_*|}.
\end{equation}

Equation \eqref{eq:full_Boltz} is typically supplemented with the inflow boundary condition: for a boundary point $\bx\in \partial\Omega$ with unit outward-pointing normal $\bn(\bx)$ and boundary velocity ${\mathbf u}_w(\bx)$, 
\begin{equation}\label{eq:bdry}
f(\bx,\bv)=f_{\text{bdry}}(\bx,\bv), \quad (\bv-{\mathbf{u}}_w(\bx))\cdot \bn(\bx)<0, \quad \bx\in \partial \Omega,
\end{equation}
where $f_{\text{bdry}}$ is either prescribed or determined by the outflow trace of the solution.
Given $f$, one can obtain the macroscopic quantities such as density, bulk velocity, and temperature via its moments:
\begin{equation}
\rho(f)(\bx)=\int_{\mathbb{R}^d}f\rd{\bv}, \quad {\mathbf u}(f)(\bx)=\frac{1}{\rho(f)}\int_{\mathbb{R}^d}f\bv\rd{\bv},\quad T(f)(\bx)=\frac{1}{d \rho(f)} \int_{\mathbb{R}^d}f|\bv-{\mathbf u}(f)|^2\rd{\bv}.
\end{equation}

In practice, the specific form of the collision kernel in \eqref{eq:kernel} cannot be determined exactly, and parametric kernels are often employed, with the parameter values determined empirically to match the correct transport properties 
observed in experiments \cite{Bird}. Indeed, Equation \eqref{eq:full_Boltz}, together with \eqref{eq:bdry}, describes a boundary-driven flow, and it is well-known that different collision kernels can produce different viscosity and diffusion coefficients (see, for example, \cite{JAH19, JAH19_1}).

This leads to an interesting inverse problem: given observational data, such as the macroscopic quantities $\rho$, $\mathbf u$, or $T$, can we determine the parameters in the corresponding collision kernel? 
To this end, we consider a kernel of the form
\begin{equation}
\label{eq:generalcollision}
B(|\bv-\bv_*|,\cos \theta;\bmu),
\end{equation}
parametrized by $\bmu\in\mM\subset\mathbb{R}^D$, where $\mM$ is a compact set.
Accordingly, the collision operator \eqref{eq:full_collision} can be written as 
\begin{align}\label{eq:Q-def}
Q(g,f;\bmu)=Q^+(g,f;\bmu)(\bv)-f(\bv)Q^-(g;\bmu)(\bv),
\end{align}
where $Q^+$ and $Q^-$ are the gain and loss terms given by
\begin{align}
Q^+(g,f;\bmu)&:=\int_{\mathbb{R}^d}\int_{\mathbb{S}^{d-1}}B(|\bv-\bv_*|,\cos \theta; \bmu)g(\bv_*')f(\bv')\rd{\sigma}\rd{\bv_*}, \\
Q^-(g;\bmu)&:=\int_{\mathbb{R}^d}\int_{\mathbb{S}^{d-1}}B(|\bv-\bv_*|,\cos \theta; \bmu) g(\bv_*)\rd{\sigma}\rd{\bv_*} .
\end{align}
Here, we express $Q$ in its bilinear form (which will be used in the following discussion) and suppress
the dependence on $\bx$. 

We then consider the following inverse problem
\begin{equation}
\label{eq:full_inverse}
\begin{aligned}
\bmu^* =\;\;  & \underset{\bmu}{\text{arg\,min}} & &\mathcal{J} (\rho(f), {\mathbf u}(f), T(f) ,\rho^{\text{data}}, {\mathbf u}^{\text{data}},  T^{\text{data}}),\\
& \text{subject~to} & & f \;\text{solves }\bv\cdot\nabla_{\bx}f=Q(f,f;\bmu) \text{ with boundary condition~\eqref{eq:bdry}},
\end{aligned}
\end{equation}
where $\mathcal{J}$ is an objective function that measures the mismatch between the predicted and observed data. The goal is to infer the parameter $\bmu$.

\begin{remark}
We note that there have been a few recent theoretical works on recovering the collision kernel in the nonlinear Boltzmann equation (e.g., \cite{LUY21,LO23,LY24}). These analytical studies establish reconstruction or uniqueness of the collision kernel from an idealized incoming-to-outgoing boundary map. Such data consist of velocity-resolved kinetic measurements under a sufficiently rich family of controlled inflow conditions. In many physical applications, however, the accessible observations are instead low-order macroscopic quantities, such as temperature, velocity, stress, or heat flux, measured under standard boundary-driven experiments. Our objective is, therefore, not to recover an unrestricted collision kernel, but to calibrate parameters in a physically motivated kernel family based on 
macroscopic observations.
\end{remark}

\section{Numerical scheme for the Boltzmann equation}\label{sec:FOM}

In this section, we describe in detail the numerical scheme (i.e., FOM) used to solve the Boltzmann equation~\eqref{eq:full_Boltz}. For simplicity, we employ a simple first-order finite difference iterative scheme coupled with the fast Fourier spectral method \cite{GHHH17} for evaluating the collision operator. 

For the present study, we consider the case  
$d=2$ and assume that the solution is homogeneous in the second spatial dimension. Consequently, $f=f(x,\bv)$, where $x\in[x_L,x_R]$ and $\bv=(v^1,v^2)\in \mathbb{R}^2$. Then Equation \eqref{eq:full_Boltz} becomes 
\begin{equation}\label{eq:Boltzmann}
v^1\partial_x f(x,\bv)=Q(f,f; \bmu)(x,\bv),
\end{equation}
where the collision operator takes the form \eqref{eq:Q-def}. 

For the boundary condition, we consider the Maxwell diffusive boundary (a special case of \eqref{eq:bdry}). For the static walls at the left boundary $x_L$ and right boundary $x_R$, with temperatures $T_L$ and $T_R$, respectively, the boundary condition is given by 
\begin{align}
\label{eq:bdry_left}
\text{when } v^1>0, \quad f(x_L,\bv)&=\rho_L\exp\left(-\frac{|\bv|^2}{2T_L} \right), \ \text{ with } \ \rho_L=\frac{\int_{v^1<0}f(x_L,\bv)(-v^1)\rd{\bv}}{\int_{v^1>0}\exp\left( -\frac{|\bv|^2}{2T_L}\right)v^1\rd{\bv}},\\
\label{eq:bdry_right}
\text{when } v^1<0, \quad
f(x_R,\bv)&=\rho_R\exp\left(-\frac{|\bv|^2}{2T_R} \right), \ \text{ with } \ \rho_R=\frac{\int_{v^1>0}f(x_R,\bv)v^1\rd{\bv}}{\int_{v^1<0}\exp\left( -\frac{|\bv|^2}{2T_R}\right)(-v^1)\rd{\bv}},
\end{align}
where $\rho_L$ and $\rho_R$ are defined to ensure zero mass flux at the boundary.

For the collision kernel, we focus on the following angularly independent form:
\begin{equation}\label{eq:B-def}
B(|\bv-\bv_*|,\cos \theta;\bmu)=\frac{b}{2\pi} |\bv-\bv_*|^a=:\Bbar(|\bv-\bv_*|;\bmu), \quad \bmu:=(a,b), \quad D=2,
\end{equation}
with the parameter range as $a\in[-1,1]$, $b\in[1,50]$, which spans regimes from hard ($a>0$) to soft potentials ($a<0$), and from kinetic/transport-dominated ($b=1$) to near-fluid/collision-dominated regimes ($b=50)$. 

If we consider the time-dependent version of Equation~\eqref{eq:Boltzmann}, $\partial_t f+v^1\partial_x f(x,\bv)=Q(f,f; \bmu)(x,\bv)$, and integrate over $x$ and $\bv$, using the mass conservation property of the collision operator together with the zero-mass-flux property at the boundary described above, one obtains 
\begin{equation}
m_0 := \int_{x_L}
^{x_R}\int_{\mathbb{R}^2}f\rd{\bv}\rd{x}=\text{constant}.
\end{equation}
This means that the total mass of the system is conserved.

Next, we describe a first-order iterative scheme for solving Equation~\eqref{eq:Boltzmann}. To set up the stage, we start with some notation. Assume that the domain $\Omega=[x_L,x_R]$ in $x$ 
is divided into $N_x$ uniform cells with $\Delta x=\frac{x_R-x_L}{N_x}$. The grid points are chosen as the midpoints of cells, that is, $x_i=x_L+(i-\frac{1}{2})\Delta x$, $\forall i\in [N_x]$.
Here $[n]=\{1, 2, \dots n\}$.
 The domain in $\bv$ is truncated to $[-L_v,L_v]^2$ with grid points $v^{1}_{j}=-L_v+(j-\frac{1}{2})\Delta v$, $v^2_{k}=-L_v+(k-\frac{1}{2})\Delta v$, $\forall j,k\in [N_v]$,
 and $\Delta v=\frac{2L_v}{N_v}$.
Associated with the mesh defined above, we introduce the discrete space
\[
\mathV = \{g: \;  g = g_{i,j,k} \ \text{at} \ (x_i, v^1_j, v^2_k), \ \forall i\in [N_x], \ \forall j,k\in [N_v]\},
\]
with the flattened representation of its function $g$  as $\gflat := \bbg(:)$, i.e., a column vector of size $\Nt=N_xN_v^2$. 
Related, we define a discrete integral and a discrete $L^2$ norm: 
\begin{equation}
    \langle h^2 \rangle_{x,\bv}  := \sum_{i,j,k} h_{i,j,k}^2\Delta x (\Delta v)^2, \quad  \|h\|_{x,\bv} := \sqrt{\langle h^2 \rangle_{x,\bv}}, \quad \forall h\in \mathV. \label{eq:def-norm}
\end{equation}
Similarly, for a grid-function $h$ only depending on $x$, we write
    $\langle h^2 \rangle_{x}  := \sum_{i} h_{i}^2\Delta x$ and $\|h\|_{x} := \sqrt{\langle h^2 \rangle_{x}}. $

Our numerical solution will be sought from the grid-function space $\mathV$. With abuse of notation, it is still represented as $\bbf$. 
For the evaluation of the collision operator, we employ the fast Fourier spectral method in \cite{GHHH17}. This method can be applied to general collision kernels of the form \eqref{eq:generalcollision}, and is therefore well suited for our purpose. In practice, it can be used as a black-box solver: given an input grid-function $\bbf\in\mathV$, 
it outputs a grid-function $\mathcal{Q}(\bbf,\bbf;\bmu) \in \mathV$, 
and similarly for $\mathcal{Q}^{\pm}$ (note that we use the notation $\mathcal{Q}$ to denote the discretized collision operator). For the grid setup discussed above, a single evaluation of the collision operator at a spatial point $x_i$ requires computational complexity $O(M N_v^3 \log N_v)$, where $M\ll N_v$ is the number of quadrature points used on the circle $\mathbb{S}^1$. Although this method is faster than the direct method\footnote{Here, the direct method refers to the direct Fourier spectral method without any acceleration strategy \cite{PR00}, whose computational complexity is $O(N_v^4)$.} for approximating the collision operator, it still represents a computational bottleneck in the overall algorithm. 

We are ready to present the iterative scheme: starting with a non-negative initial guess $\bbf^0\in\mathV$ and $\langle \bbf^0 \rangle_{x,\bv}=m_0>0$ (e.g.,  a spatially uniform Maxwellian),  we proceed as follows to update the $(n+1)$-th iterate $\bbf^{n+1}\in\mathV, n\geq 0$:
\begin{itemize}
\item 
For any $j\in [N_v]$ such that $v^1_j\geq 0$ and $\forall k\in [N_v]$, we sweep from $i=1$ to $i=N_x$ and compute $f_{i,j,k}^{n+1}$ based on 
\begin{equation}
\label{eq:sweep+}
v^1_j\frac{f_{i,j,k}^{n+1}-f_{i-1,j,k}^{n+1}}{\Delta x}=\mathcal{Q}^+(\bbf^n,\bbf^n;\bmu)_{i,j,k}-f_{i,j,k}^{n+1}\mathcal{Q}^-(\bbf^n;\bmu)_{i,j,k}. 
\end{equation}
Here, the ghost value $f_{0,j,k}^{n+1}:=\mL_{bdry}^L(f^n; T_L, v_j^1, v_k^2)$ imposes the boundary condition~\eqref{eq:bdry_left}, with $\mL_{bdry}^L$ defined as follows with $g\in \mathV$,
\begin{equation}
\mL_{bdry}^L(g; T_L, v_j^1, v_k^2)
=\bar{\rho}_L\exp\left(-\frac{(v^1_j)^2+(v^2_k)^2}{2T_L} \right), \;\;  \bar{\rho}_L=\frac{\sum_{v^1_l<0}g_{1,l,k}(-v^1_l)}{\sum_{v^1_l>0} \exp\left( -\frac{(v^1_l)^2+(v^2_k)^2}{2T_L}\right)v^1_l}.
\end{equation}

\item For any $j\in [N_v]$ such that $v^1_j<0$ and $\forall k\in [N_v]$, we sweep from $i=N_x$ to $i=1$ and compute $f_{i,j,k}^{n+1}$ based on  
\begin{equation}
\label{eq:sweep-}
v_j^1\frac{f_{i+1,j,k}^{n+1}-f_{i,j,k}^{n+1}}{\Delta x}=\mathcal{Q}^+(\bbf^n,\bbf^n;\bmu)_{i,j,k}-f_{i,j,k}^{n+1}\mathcal{Q}^-(\bbf^n;\bmu)_{i,j,k}.
\end{equation}
Here, the ghost value $f_{N_x+1,j,k}^{n+1}:=\mL_{bdry}^R(f^n; T_R, v_j^1, v_k^2)$ imposes the boundary condition \eqref{eq:bdry_right}, with $\mL_{bdry}^R$ defined as follows with $g\in \mathV$,
\begin{equation}
\mL_{bdry}^R(g; T_R, v_j^1, v_k^2)
=\bar{\rho}_R\exp\left(-\frac{(v^1_j)^2+(v^2_k)^2}{2T_R} \right), \;\; \bar{\rho}_R=\frac{\sum_{v^1_l>0}g_{N_x,l,k}v^1_l}{\sum_{v^1_l<0} \exp\left( -\frac{(v^1_l)^2+(v^2_k)^2}{2T_R}\right)(-v^1_l)}.
\end{equation}

\item A linear scaling is applied to preserve the mass: $\bbf^{n+1}\leftarrow  \frac{\bbf^{n+1}}{\langle \bbf^{n+1} \rangle_{x,\bv}}m_0$.
\end{itemize}

The iteration is terminated when the discrete $L^2$ norm of the residual of $\bbf^{n+1}$,  $\|\Res(\bbf^{n+1}; \bmu)\|_{x,\bv}$,  falls below a prescribed threshold. Here, the residual of any $g\in \mathV$ (with respective to a given  $\bmu$) at $(x_i, v^1_j, v^2_k)$ is defined as 
\begin{equation}
\label{eq:res}
\Res(\bbg; \bmu)_{i,j,k}
:=\mL(g)_{i,j,k}-\mathcal{Q}(\bbg,\bbg;\bmu)_{i,j,k},
\end{equation}
where $\mL: \mathV\mapsto \mathV$ denotes the linear discrete operator for the transport term $v^1\partial_x$ given by
\begin{equation}
\mL(g)_{i,j,k}:=\frac{v^1_j+|v^1_j|}{2}\frac{g_{i,j,k}-g_{i-1,j,k}}{\Delta x}+\frac{v^1_j-|v^1_j|}{2}\frac{g_{i+1,j,k}-g_{i,j,k}}{\Delta x}.
\end{equation}
Again, we set 
   $ g_{0,j,k}=\mL_{bdry}^L(g; T_L, v_j^1, v_k^2)$ when $v_j^1\geq 0$ and $g_{N_x+1,j,k}=\mL_{bdry}^R(g; T_R, v_j^1, v_k^2)$ when $v_j^1< 0$.
As in the continuous case, the discrete boundary operators $\mL_{bdry}^L$ and $\mL_{bdry}^R$ remain linear in $g$. This is crucial when we consider the proposed ROMs, which are defined on linear reduced spaces.

\begin{remark} It is worth pointing out that, with the upwind discretization for the spatial derivative and the sweeping strategies in our iterative scheme,  \eqref{eq:sweep+} is equivalent to a simple update 
\begin{equation}
f_{i,j,k}^{n+1} = \frac{v^1_jf_{i-1,j,k}^{n+1}+\Delta x \mathcal{Q}^+(\bbf^n,\bbf^n;\bmu)_{i,j,k}}{v^1_j+\Delta x\mathcal{Q}^-(\bbf^n;\bmu)_{i,j,k}},
\end{equation}
while  \eqref{eq:sweep-} is
\begin{equation}
f_{i,j,k}^{n+1} = \frac{-v^1_jf_{i+1,j,k}^{n+1}+\Delta x \mathcal{Q}^+(\bbf^n,\bbf^n;\bmu)_{i,j,k}}{-v^1_j+\Delta x\mathcal{Q}^-(\bbf^n;\bmu)_{i,j,k}}.
\end{equation} 
\end{remark}

\begin{remark}
The scheme presented in this section is also referred to as source iteration, following the terminology used in solving linear transport equations. It is well-known that this method suffers from slow convergence in the fluid (collision-dominated) regime. In such cases, the preconditioning technique proposed in \cite{CDH25} can be used to accelerate convergence. Furthermore, when the dimension $d$ is high, the full grid-based method can become computationally expensive, making the dynamical low-rank method particularly advantageous \cite{HW22}. Since the FOM is not the primary focus of this work, we do not employ these techniques.
\end{remark}

The numerical solution by the iterative scheme will be referred to as the FOM solution. As the high-fidelity solution, it will be regarded as the ground truth hereafter and used to build our ROM in the next section. %
To highlight the dependence on the parameter, we also write it as $\bbf(\cdot\,; \bmu)=\text{FOM}(\bmu)$, while its flattened is denoted as $\fflat_{\bmu}$.

\section{Reduced order model  for the Boltzmann equation}\label{sec:ROM}

Our ROM is constructed following the reduced basis method (RBM) framework \cite{hesthaven2016certified, haasdonk2017reduced}, which consists of an offline stage and an online stage. In the offline stage, a nested sequence of low-dimensional approximation spaces $\mathV_1\subset \mathV_2\subset \cdots \subset \mathV_N$ is built iteratively, with each termed as an RB space and given as $\mathV_n=\text{span}\{\fFOM(\cdot\,;\bmu_1),\fFOM(\cdot\,;\bmu_2),\dots, \fFOM(\cdot\,;\bmu_n)\}$, where the parameter values $\{\bmu_n\}_{n=1}^N$ are selected in a greedy fashion; in the online stage, reduced solutions at many parameter samples are sought from the terminal RB space $\mathV_N$. Given our parametric PDE model is nonlinear, special care is needed to ensure the overall cost efficiency both online and offline. We particularly utilize the quadratic nature of the collision operator in conjunction with a separable approximation for the collision kernel. Furthermore, we incorporate some fundamental physical properties of the solutions.

\subsection{Efficient evaluation of the parametric collision operator}
\label{sec:EIM}

The overall efficiency of our reduced surrogate solver and its application to the inverse problem is greatly facilitated by using a
separable approximation of the parametric collision kernel 
$B(\cdot\,,\cdot\,; \bmu)$ along with the precomputation of some collision terms.  
Note that the collision kernel \eqref{eq:B-def} considered in this work is $\Bbar(r; \bmu)$ with $r=|\bv-\bv_*|$ and $\bmu=(a,b)$. 
In addition, 
\begin{equation}\label{eq:B-kappa}
\Bbar(r; \bmu)|_{\bmu=(a,b)}=b\underbrace{\Bbar(r; \mueim)|_{{\mueim}=(a,1)}}_{\kappa(r; a)}, 
\text{ with } \kappa(r; a) = \frac{1}{2\pi}r^a.
\end{equation}

Note that $b$ in \eqref{eq:B-kappa} appears only as a multiplicative factor. Thus, the task of finding a separable approximation for the collision kernel $\Bbar(r; \bmu)$ reduces to approximating 
$\kappa(r;a)$ by 
a separable function $\kappa^S(r; a)$, namely,
\begin{equation}
    \kappa(r; a)\approx \kappa^S(r; a):=\sum_{m=1}^S \beta_m(a) \kappa(r; a_{m}),
    \label{eq:sep:approx}
\end{equation}
and equivalently,
\begin{equation}
    \Bbar(r; \bmu)|_{\bmu=(a,b)}\approx %
    b\sum_{m=1}^S \beta_m(a) \Bbar(r; \mueim_m)|_{{\mueim_m}=(a_m,1)},
    \label{eq:sep:approx1}
\end{equation}
with a moderate $S$. This is achieved by the empirical interpolation method (EIM)  \cite{barrault2004empirical, hesthaven2016certified}. Specifically, a set of parameter points  $\{a_m\}_{m=1}^S\subset [a_-, a_+]$ and a set of interpolation points $\{r_m\}_{m=1}^S\subset [r_-, r_+]$ are selected by a greedy algorithm (until certain stopping criterion, e.g., based on  interpolation errors, is satisfied), where $a_{k+1}$ (with $k<S$) is chosen as the next parameter point for which  $\kappa(r; a_{k+1})$ is the worst approximated by its interpolation\footnote{The interpolation error is measured in $\ell^p$-norm in general, with $p=\infty$ in this work.} in $\text{span}\{\kappa(r; a_m)\}_{m=1}^k$ based on the interpolation points  $\{r_m\}_{m=1}^k$, while $r_{k+1}$ is chosen as the point where the interpolation error of  $\kappa(r; a_{k+1})$  is maximal. Here $[a_-, a_+]=[-1,1]$ and $[r_-, r_+]=[0,R]$ (with $R=9$) are the relevant ranges of $a$ and $r$ used in numerical examples. 
  
 Once $\{a_m\}_{m=1}^S$ and $\{r_m\}_{m=1}^S$ are available, we compute and store the interpolation matrix $\bG\in \mathbb R^{S\times S}$,
 with $G_{ij}=\kappa(r_i; a_j)$, along with its inverse $\bG^{-1}$. (The invertibility of $\bG$ is guaranteed by the stopping criterion.)  For any given $a\in[a_-, a_+]$, the coefficient vector $\bbeta(a)=\big(\beta_m(a)\big)\in \mathbb R^S$ in \eqref{eq:sep:approx} is determined via interpolation, namely, $\kappa^S(r_j;a)=\kappa(r_j;a), j\in[S]$.
 Algebraically, this is equivalent to 
\begin{equation}\label{eq:beta-formula}
    \bbeta(a)=\bG^{-1} \; [\kappa(r_1;a), \kappa(r_2;a), \cdots, \kappa(r_S;a)]^T.
\end{equation}
 
 Standard implementation of EIM can be referred to, e.g., in \cite{hesthaven2016certified}. Below are our two specific considerations.  
\begin{itemize}
    \item \textbf{Shifted/transformed $\kappa(r; a)$.}  In our work, the EIM is applied to $r\kappa(r;a)=\kappa(r;a+1)$ instead of $\kappa(r,a)$.
    Unlike in  \cite{antil2019reduced}, a small neighborhood of $r=0$ does not need to be excluded. Although $\kappa(r;a)$ alone displays a singularity at $r=0$ when $a\in[-1,0)$,  we indeed work with $r\kappa(r;a)$ in the collision solver (we refer to \cite{GHHH17} for more details),  hence only require  $r\kappa^S(r,a)$ to be a sufficiently accurate approximation for $r\kappa(r;a)$.
\item \textbf{Snapshot EIM basis.} Separable approximations by EIM are usually expressed in a different type of basis \cite{hesthaven2016certified}, as its interpolation matrix is better conditioned than $\bG$. Given that $\kappa^S(r; a)$ is used together with a black-box collision solver, it is crucial for us to work with the snapshot type EIM basis $\{\kappa(r; a_{m})\}_{m=1}^S$ to represent  $\kappa^S(r; a)$ as in \eqref{eq:sep:approx}.   Also see \Cref{sec:num} for more discussions.
\end{itemize}

Next, we illustrate how the computation of parametric collision operator can benefit from the separable kernel approximation. 
In our application, we need to evaluate $\mathcal{Q}(h,h; \bmu)$,  where $h =\sum_{k=1}^n c_k h_k$, and  $\{h_k\}_{k=1}^n$ is a given set of grid functions in $\mathV$. We are particularly interested in such evaluations for a large number of coefficient vectors $\bc\in \mathbb{R}^n$ and parameter values $\bmu=(a,b)$.  For each relevant $a$, assume $\bbeta(a)$ in \eqref{eq:beta-formula} is available. By using the intrinsic
quadratic structure of the collision operator \eqref{eq:Q-def} and the separable approximation of the kernel,  we have
\begin{align}
\left.\mathcal{Q}\left(h, h; \bmu\right)\right|_{\bmu=(a,b)}  = & \, \mathcal{Q}\left(\sum_{i=1}^n c_i h_i, \sum_{i=1}^n c_i h_i; \bmu\right)  = \sum_{k,l=1}^n c_k c_l \mathcal{Q}(h_k,h_l; \bmu) \notag \\
\mathring{\approx}  & \,    
b\sum_{k,l=1}^n c_k c_l  \left[\sum_{m=1}^S \beta_m(a) \underbrace{\mathcal{Q}(h_k,h_l; \mueim_m)}_{\text{precomputed}}\big|_{\mueim_m=(a_m,1)} \right], \label{eq:Q-eim}
\end{align}
where the notation $\mathring{\approx}$, used here and below, indicates that the approximation on the right-hand side is used in our actual algorithm.
This implies  that once the collisions terms $\{\mathcal{Q}(h_k,h_l; \mueim_m)\big|_{\mueim_m=(a_m,1)} \}_{m=1}^S$ are precomputed, we can evaluate $\mathcal{Q}(h, h; \bmu)$ efficiently for many instances of $\bc$ and $\bmu$.

\subsection{Reduced order model: offline stage}\label{sec:offline}

To build the low-dimensional RB spaces to approximate the solutions to the Boltzmann equation at many parameter samples, namely,
$\{ \fFOM(\cdot\,; \bmu) : \bmu \in \Mrange \subseteq \mM \}$, we follow a greedy approach that starts  with a training parameter set $\Mtrain \subset \mM$. 
 Given the current RB space ${\mathbb V}_n=\textrm{span}\{\fFOM(\cdot\,;\bmu_i)\}_{i=1}^n$, with $\bmu_i\in \Mtrain, \ i\in[n]$, the greedy procedure seeks the parameter $\bmu_{n+1}\in\Mtrain$ whose solution  
is worst represented by the current ROM basis. More specifically, for each $\bmu\in\Mtrain$, we consider its ``best linear representation’’ in ${\mathbb V}_n$, referred to as the reduced solution $\ffROM^n(\cdot\,;\bmu)=\text{ROM}(\mathV_n;\bmu)$ and to be defined in \eqref{eq:rom:ansatz}–\eqref{eq:c-constrained-opt}. The parameter $\bmu_{n+1}$ is then selected from the remaining training samples as the one for which the residual norm \eqref{eq:res} of this reduced solution is the largest, and $\fFOM(\cdot\,;\bmu_{n+1})$ is subsequently added to the RB space.

 Unlike the widely used POD approach to build low-dimensional approximation spaces in the MOR framework, with the greedy strategy, the minimum number of FOM solutions needs to be computed  during our offline stage. This is relevant to the overall efficiency of our MOR strategy for kinetic models, as the offline cost is dominated by that of the FOM solves, especially near the fluid regime; also see \Cref{sec:num-rom}.

Let $\Mtrain\subset\mM$ be a sufficiently large training parameter set, $N_{\text{max}}$ be the maximal number of greedy iterations, $\tau_{tol}>0$ be a prescribed tolerance;  Set  $\mathV_0=\emptyset$ and  $\rbF_0=\emptyset$.   The greedy algorithm of the offline stage of the ROM is summarized below: 

\medskip
\noindent\textbf{Prestage of EIM selection:}
Following the procedure in \Cref{sec:EIM}, we use EIM to select parameter $a$'s $\{ a_m \}_{m=1}^S$ and the interpolation points $\{ r_m \}_{m=1}^S$, formulate the corresponding EIM parameters $\{\mueim_m :=(a_m, 1)\}_{m=1}^S$, 
and compute the inverse interpolation matrix $\bG^{-1}$. 
Precompute $\{\bbeta(a) = (\beta_m(a)): (a,b) \in \Mtrain \}$ using \eqref{eq:beta-formula}.

\medskip
\noindent\textbf{Initialization:} The first parameter value $\bmu_1$ is chosen randomly from $\Mtrain$. Compute $\fFOM(\cdot\,;\bmu_1)=\text{FOM}(\bmu_1)$, and its total mass $\mass = 
\left\langle \fFOM(\cdot\,;\bmu_1) \right\rangle_{x,\bv}$. Set the RB space ${\mathbb V}_1=\textrm{span}\{\fFOM(\cdot\,;\bmu_1)\}$, and its matrix representation $\rbF_1=[\fflat_{\bmu_1}]\in\mathR^{\Nt}$.

\medskip
\noindent\textbf{Greedy basis generation:} 
for $n=1,2,\cdots, N_{\text{max}}$

\begin{enumerate}
    \item \textbf{Orthogonalization}: find the orthogonal basis $\{ u_i \}_{i=1}^n$ for the current RB space ${\mathbb V}_n$. This is achieved by performing the economy-size QR decomposition of its associated matrix representation $\rbF_n$, namely 
\begin{equation}\label{eq:QR}
[\bU, \bR] = \texttt{QR}(\rbF_n,  \texttt{"econ"}).
\end{equation}  
    The basis function $u_i\in \mathV_n$ is the corresponding grid-function of the $i$-th column of  $\bU$.  Here $\bU\in \mathR^{\Nt\times n}$, satisfying $\bU^T\bU=\bI_{n\times n}$, and $\bR\in \mathR^{n\times n}$ is upper triangular.  Given the nested structure $\rbF_n=[\rbF_{n-1}, \fflat_{\bmu_n}]\in\mathR^{\Nt\times n}$ 
    and with QR, one only needs to compute incrementally  the last column  of  $\bU$ corresponding to $\bu_n$.  

\item   \textbf{Precomputation related to transport and collision operators}: The computational efficiency benefits from the following two subsets of $\mathV$ being precomputed: 
\[\{\mathcal{Q}(\bu_i,\bu_j; \mueim_m):\  i,j\in[n],\ m\in[S]\},\quad\text{and}\quad\{\mL(u_i):\ i\in [n]\}.\]
By utilizing the nested structure of the RB space $\mathV_{n-1}\subset \mathV_n$, one only needs to compute incrementally the new terms, namely
$\mL(\bu_n)$, and
\begin{equation}\label{eq:Q-extra}
  \{\mathcal{Q}(\bu_i,\bu_n; \mueim_m), \mathcal{Q}(\bu_n, \bu_i; \mueim_m): i\in[n], \, m\in[S]\}.  
\end{equation}
         \item  \textbf{Find the next parameter value $\bmu_{n+1}$}, provided that the termination condition is not met, that is, the spectral ratio %
         of $\rbF_n$ is {above} the tolerance $\tau_{tol}$:\footnote{Following the work in \cite{peng2022reduced}, the spectral ratio of $\rbF_n$ is defined as %
         $\sigma_{n}/(\sum\limits_{i=1}^{n} \sigma_i)$, where $\sigma_i$ is the $i$-th largest singular value of $\rbF_n$. } %

\begin{itemize}
    \item 
         \textbf{For any $\bmu\in \Mtrain\backslash \{\bmu_i\}_{i=1}^n$, find the reduced solution from $\mathV_n$} 
         \begin{equation}
            \ffROM^n(\cdot\,; \bmu):=\sum_{i=1}^n c_i(\bmu) \bu_i\quad  (\text{also denoted as}\; \text{ROM}(\mathV_n; \bmu))
        \label{eq:rom:ansatz}
         \end{equation}  via residual minimization. This is equivalent to solving  the coefficient vector  $\bc(\bmu) = (c_i(\bmu))\in \RR^{n}$ from the following optimization problem while  preserving the total mass:  %
\begin{equation}\label{eq:c-constrained-opt}
\begin{aligned}
\bc(\bmu) :=
\underset{\bc\in \RR^n}{\text{arg\,min}} \;\;& 
\left\|\Res\left(\sum_{i=1}^n c_i \bu_i; \bmu\right)\right\|_{x,\bv}^2
\\
\text{subject to}\; & 
\left\langle \sum_{i=1}^n c_i \bu_i \right\rangle_{x,\bv} =\mass.
\end{aligned}
\end{equation}
The objective value at the optimizer is denoted as 
 \begin{equation}\label{eq:minRes}
     \MinRes(\bmu):=\left\|\Res\left(\sum_{i=1}^n c_i(\bmu) \bu_i; \bmu\right)\right\|_{x,\bv}^2.
 \end{equation}

\item \textbf{Select the next parameter $\bmu_{n+1}$ and update the RB space}.  The next parameter $\bmu_{n+1}$ is determined as the one whose reduced solution is approximated least well by the current RB space $\mathV_n$, in the sense that 
$$\bmu_{n+1}: =\underset{\bmu\in \Mtrain\backslash \{\bmu_i\}_{i=1}^n}{\text{arg\,max}} \MinRes(\bmu).
$$
Compute $\fFOM(\cdot\,;\bmu_{n+1})=\text{FOM}(\bmu_{n+1})$, update the RB space ${\mathbb V}_{n+1}={\mathbb V}_{n} \bigoplus \{\fFOM(\cdot\,;\bmu_{n+1})\}$, and its matrix representation $\rbF_{n+1}=[\rbF_{n}, \fflat_{\bmu_{n+1}}]$.  
\end{itemize}

\end{enumerate}

We now want to highlight and elaborate on some technical details and considerations related to, e.g.,  the structure-preservation and efficiency, of $\text{ROM}({\mathV}_n; \bmu)\protect$ in \eqref{eq:rom:ansatz}:
\begin{itemize} 

\item \textbf{Orthogonalization.} When working with the reduced space $\mathV_n$, or equivalently its algebraic representation $\rbF_n$, orthogonalization (e.g., via QR as in \eqref{eq:QR}) is  crucial to improve the conditioning and accuracy of the optimization problem.

\item \textbf{Conservation.}  When finding the reduced solution $\ffROM^n(\cdot\,; \bmu)$ from $\mathV_n$, the mass conservation is strictly enforced in the optimization problem \eqref{eq:c-constrained-opt}, %
along with an initial guess $\bc = \bR [1/n, \ldots, 1/n]^\top$, that corresponds to using $(\sum_{i=1}^n \fflat_{\bmu_i})/n$ to initialize $\ffROM^n(\cdot\,; \bmu)$ to solve the optimization task.
Recall that the matrix $\bR$ is from QR in \eqref{eq:QR}, connecting the orthonormal basis with the original $\fflat$ solution.
  In the case when the greedily selected snapshots  $\fflat_{\bmu_i}, i\in [n]$ are all positive, such initialization provides a positive initial guess for the reduced solution. 

\item \textbf{Efficient computation of reduced solutions via accelerated residual evaluation.} 
The reduced solution $\ffROM^n(\cdot\,; \bmu)=\text{ROM}(\mathV_n; \bmu)$ is determined in \eqref{eq:rom:ansatz}-\eqref{eq:c-constrained-opt} via solving a constrained optimization problem whose objective involves the residual of any grid-function in ${\mathbb V}_n$. This requires repeated evaluations of the objective and the residual inside throughout the optimization iterations, which are
computationally expensive due to the nonlinear collision term \eqref{eq:Q-def} in the Boltzmann equation. %
In our algorithm, we utilize the intrinsic quadratic structure of the collision operator hence the quadratic dependence of $\mathcal{Q}\left(\sum_{i}^n c_i \bu_i, \sum_{i}^n c_i \bu_i; \bmu\right)$ on the coefficient vector $\bc$, as well as the separable approximation \eqref{eq:Q-eim} of the kernel. Thus, the squared $L^2$-norm of the residual can be approximated and computed as
\begin{equation}\label{eq:res-eim}
\left\|\Res\left(\sum_{i}^n c_i \bu_i; \bmu\right) \right\|_{x,\bv}^2
\mathring{\approx}
\left\|  \sum_{i=1}^n c_i \mL(\bu_i)  - 
b\sum_{k,l=1}^n c_k c_l   \left[\sum_{m=1}^S \beta_m(a) \mathcal{Q}(\bu_k,\bu_l; \mueim_m) \right]
\right\|_{x,\bv}^2.
\end{equation}
The reduced solution can therefore be computed efficiently using the precomputed terms $\mL(u_i)$, $\mathcal{Q}(\bu_k,\bu_l; \mueim_m)$ and $\beta_m(a)$. %

\item \textbf{Constrained optimization problem.}  To solve the constrained optimization problem \eqref{eq:c-constrained-opt} numerically, we use the \texttt{MATLAB} function \texttt{fmincon} with the sequential quadratic programming (SQP) algorithm (see~\cite[Chapter 18]{nocedal1999numerical}). The gradient of the objective function with respect to $\bc$ is given 
analytically as follows
and directly used in our implementation,
again due to the quadratic relationship of $\Res\left(\sum_{i=1}^n c_i \bu_i; \bmu\right) $ in $\bc$:
\begin{equation}\label{eq:dc_Res}
\begin{aligned}
&\partial_{c_j} \left\| \Res\left(\sum_{i=1}^n c_i \bu_i; \bmu\right) \right\|_{x, \bv}^2\\
\mathring{\approx} & \partial_{c_j} 
\left\|  \sum_{i=1}^n c_i \mL(\bu_i)  - b
\sum_{k,l=1}^n c_k c_l   \left(\sum_{m=1}^S \beta_m(a) \mathcal{Q}(\bu_k,\bu_l; \mueim_m) \right)
\right\|_{x,\bv}^2
\\
= & 
2 \left\langle
\left[ \sum_{i=1}^n c_i \mL(\bu_i)  - 
b\sum_{k,l=1}^n c_k c_l \left(\sum_{m=1}^S \beta_m(a) \mathcal{Q}(\bu_k,\bu_l; \mueim_m) \right)
\right] \right. \\
&  \qquad
\left.\cdot
\left[ 
\mL( \bu_j)  -b 
\sum_{i=1}^n c_i \left( 
\sum_{m=1}^S \beta_m(a) \Big(\mathcal{Q}(\bu_j,\bu_i; \mueim_m) 
+  \mathcal{Q}(\bu_i,\bu_j; \mueim_m) \Big)
\right)
\right] \right\rangle_{x,\bv}
.
\end{aligned}
\end{equation}
\end{itemize}

Given the steps and details discussed above, the offline stage of the ROM is summarized in \Cref{alg:offline}.

\begin{algorithm}[h!]
\caption{Offline stage of ROM for
Boltzmann equation, using precomputed collision terms and EIM kernel approximation}\label{alg:offline}
\begin{algorithmic}[1]

\STATE \textbf{Input:} Training parameter set $\Mtrain$, maximum number of iterations $N_{\max}$, spectral tolerance $\tau_{\mathrm{tol}}$

\STATE \textbf{Prestage (EIM):} 
Select EIM parameter %
$\Meim=\{ \mueim_m = (a_m,1) \}_{m=1}^S$, compute inverse interpolation matrix $\bG^{-1}\in\mathbb{R}^{S\times S}$, and precompute $\{\bbeta(a) = (\beta_m(a))\in \mathbb{R}^{S}: (a,b) \in \Mtrain \}$ using \eqref{eq:beta-formula}.

\STATE \textbf{Initialization:}
Choose randomly $\bmu_1\in\Mtrain$, compute $\fFOM(\cdot\,;\bmu_1)=\text{FOM}(\bmu_1)$ and $\mass = 
\left\langle \fFOM(\cdot\,;\bmu_1) \right\rangle_{x,\bv}$. Set  ${\mathbb V}_1=\textrm{span}\{\fFOM(\cdot\,;\bmu_1)\}$ and  $\rbF_1=[\fflat_{\bmu_1}]$.

\FOR{$n=1,\dots,N_{\max}$}

    \STATE \textbf{Orthogonalization:}
    Compute an orthonormal basis $\{ u_i \}_{i=1}^n$ 
    of $\mathV_n$ via (e.g., incremental implementation as in \cite{daniel1976reorthogonalization}) of QR,  $[\bU, \bR] = \texttt{QR}(\rbF_n,\texttt{"econ"})$ 
    \STATE \textbf{Precomputation:}
    Update precomputed transport and collision terms incrementally with new terms $\mL(u_n)$ and $\{\mathcal{Q}(\bu_i,\bu_n; \mueim_m), \mathcal{Q}(\bu_n, \bu_i; \mueim_m): i\in[n],
    \, m\in[S]
    \}$

    \STATE \textbf{Stopping criterion:}
    Compute spectral ratio $r_{n}$ of $\rbF_{n}$; \textbf{if} $r_{n} \le \tau_{\mathrm{tol}}$ \textbf{stop}

    \STATE \textbf{Greedy selection:}
\FOR{each $\bmu \in \Mtrain \setminus \{\bmu_i\}_{i=1}^n$}
    \STATE Compute ROM solution $\ffROM^n(\cdot\,;\bmu)=\sum_{i=1}^n c_i(\bmu) \bu_i$ by solving the optimization problem \eqref{eq:c-constrained-opt} to find coefficients $\bc(\bmu)$ minimizing the norm of the residual (computed using the precomputed terms and EIM  via \eqref{eq:res-eim})
    \STATE Evaluate the objective $\MinRes(\bmu)$ via \eqref{eq:minRes}
\ENDFOR
\STATE Choose next parameter $\bmu_{n+1} = \arg\max_{\bmu} \MinRes(\bmu)$,
    compute $\fFOM(\cdot\,;\bmu_{n+1})=\text{FOM}(\bmu_{n+1})$, update ${\mathbb V}_{n+1}={\mathbb V}_{n} \bigoplus \{\fFOM(\cdot\,;\bmu_{n+1})\}$ and $\rbF_{n+1}=[\rbF_{n}, \fflat_{\bmu_{n+1}}]$

\ENDFOR

\STATE \textbf{Output:} Reduced dimension $N=n$, RB basis $ \{u_i\}_{i=1}^N$, precomputed terms $\{ \mL (u_i):\ i\in[N] \}$ and $\{\mathcal{Q}(\bu_i,\bu_j; \mueim_m):\ i,j\in[N], \ m\in[S]
\}$, EIM points $\{\mueim_m,r_m\}_{m=1}^S$ and inverse interpolation matrix $\bG^{-1}$

\end{algorithmic}
\end{algorithm}

\subsection{Reduced order method: online stage}\label{sec:online}

Let $\mathV_N$ be the terminal RB space from the offline stage of the reduced order method, with the orthonormal basis $\{ u_i\}_{i=1}^N$. Upon termination,  we have the precomputed  quantities $\{\mL(\bu_k): k\in[N]\}$ and $\{\mathcal{Q}(\bu_k,\bu_l; \mueim_m): k,l\in[N], m\in[S]\}$.  Let $\Mrange\subset\mM$ be the set of parameter values under consideration. 
For any $\bmu=(a,b)\in\Mrange$, the reduced solution $\ffROM(\cdot\,; \bmu)=\sum_{i=1}^N c_i(\bmu) \bu_i=\text{ROM}(\mathV_N; \bmu)$ can be found via the following two steps: 
\begin{enumerate}
    \item Step 1: Compute $\bbeta(a)$ using \eqref{eq:beta-formula}.
    \item Step 2: the coefficient vector $\bc(\bmu)$ of $\ffROM(\cdot\,; \bmu)$ is determined by minimizing the residual, computed through the optimization \eqref{eq:c-constrained-opt} with the greatly simplified residual in \eqref{eq:res-eim}.
\end{enumerate}

The full algorithm for online stage of ROM is summarized in \Cref{alg:online}.

\begin{algorithm}[h]
\caption{Online stage of ROM to compute 
the RB solution in ${\mathbb V}_{N}$ to Boltzmann equation, 
using precomputed collision terms and EIM kernel approximation}
\begin{algorithmic}[1]\label{alg:online}

\STATE \textbf{Input:} parameter $\bmu=(a,b)$, RB basis $\{u_i\}_{i=1}^N$, precomputed  terms $\{ \mL (u_i):\ i\in [N]\}$
and $\{\mathcal{Q}(\bu_i,\bu_j; \mueim_m):\  i,j\in [N],\ m\in[S]
\}$, EIM points $\{\mueim_m,r_m\}_{m=1}^S$ and inverse interpolation matrix $\bG^{-1}$

\STATE \textbf{EIM evaluation:} Evaluate $\{\kappa(r_m; a)\}_{m=1}^S$ and compute $\bbeta(a) = (\beta_m(a))$ using \eqref{eq:beta-formula} 

    \STATE \textbf{Optimization:} Compute ROM solution $\ffROM(\cdot\,;\bmu)=\sum_{i=1}^N c_i(\bmu) \bu_i$ by solving the optimization problem \eqref{eq:c-constrained-opt} to find coefficients $\bc(\bmu)$ minimizing the norm of the residual (computed using the precomputed terms and EIM  via \eqref{eq:res-eim})

\STATE \textbf{Output:} ROM solution $\ffROM(\cdot\,;\bmu)=\sum_{i=1}^N c_i(\bmu) \bu_i$ and the corresponding coefficient vector $\bc(\bmu)$
\end{algorithmic}
\end{algorithm}

\section{Solving inverse problem using RB space and ROM}
\label{sec:IVP}

Building on the reduced model, we now consider the inverse problem \eqref{eq:full_inverse} to reconstruct the parameters in the collision kernel from macroscopic observations. In particular, we focus on a thermally driven flow (Fourier flow), in which  the temperatures $T_L$ and $T_R$ are prescribed at the two boundaries $x_L$ and $x_R$, with $T_L\ne T_R$. Different choices of the collision parameters $\bmu=(a,b)$ give rise to different temperature profiles in the interior of the domain. 

Assume that temperature observation data, denoted by $T^{\mathrm{data}}(x)$, are available throughout the physical domain $\Omega=[x_L,x_R]$. The goal is to reconstruct the parameter $\bmu$ in the collision kernel. This can be formulated as a least-squares optimization problem by minimizing the data misfit between the observed temperature and the temperature predicted by the governing equation:
\begin{equation}\label{eq:IP-Boltz}
\begin{aligned}
\bmu^* =\;\;  & \underset{\bmu=(a,b)}{\text{arg\,min}} & &\| T(f) - T^{\text{data}} \|^2_{L^2}\\
& \text{subject to} & & v^1\partial_x f=Q(f,f; \bmu) \text{\,\,and boundary condition~ \eqref{eq:bdry_left}-\eqref{eq:bdry_right}}.
\end{aligned}
\end{equation}

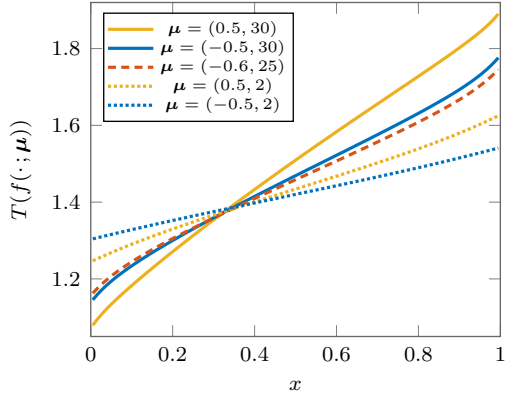
\begin{figure}[h]
    \centering
    \begin{tikzpicture}
    \begin{axis}[
    width=7cm,
    height=6cm,
    xmin=0, xmax=1,
    ymin=1.05, ymax=1.92,
    xlabel={$x$},
    ylabel={$T(f(\cdot\,;\bmu))$},
    tick label style={font=\scriptsize},
    label style={font=\scriptsize},
    legend style={
        fill=white,
        fill opacity=0.50,
        text opacity=1,
        font=\tiny,
        at={(0.03,0.97)},
        anchor=north west,
        row sep=-3pt,
        inner sep=1pt
    },
    grid style={gray!25},
    major grid style={gray!35},
    line width=0.5pt,
    axis line style={black!60},
]
    \addplot[matlabyellow, very thick] table[x=x,y=T_a_p05_b30] {data/temperature-observation-boltz-clean.txt};
\addlegendentry{$\bmu=(0.5,30)$}
    \addplot[matlabblue, very thick] table[x=x,y=T_a_m05_b30] {data/temperature-observation-boltz-clean.txt};
    \addlegendentry{$\bmu= (-0.5,30)$}
    \addplot[matlaborange, very thick, densely dashed] table[x=x,y=T_a_m06_b25] {data/temperature-observation-boltz-clean.txt};
    \addlegendentry{$\bmu=(-0.6,25)$}
	    \addplot[matlabyellow, very thick, densely dotted] table[x=x,y=T_a_p05_b2] {data/temperature-observation-boltz-clean.txt};
	    \addlegendentry{$\bmu=(0.5,2)$}   
         \addplot[matlabblue, very thick, densely dotted] table[x=x,y=T_a_m05_b2] {data/temperature-observation-boltz-clean.txt};
	    \addlegendentry{$\bmu=(-0.5,2)$}
	\end{axis}%
     \end{tikzpicture}
    \caption{Temperature profile for different parameter $\bmu$.}
    \label{fig:T}
\end{figure}

Since, to our knowledge, there are no theoretical results for this inverse problem,
we first perform a parameter study to assess the difficulty of the reconstruction problem. We assume $[x_L,x_R]=[0,1]$, $T_L=1$ and $T_R=2$, and plot the temperature profiles $T(f(\cdot\,;\bmu))$ obtained from the FOM solver for different values of $\bmu$ in \Cref{fig:T}. It can be observed that the temperature profiles differ only subtly when the parameter changes from $\bmu=(-0.5,30)$ to $\bmu=(-0.6,25)$. This illustrates the difficulty of the inverse problem: a small perturbation in the observation data may correspond to a large change in the underlying parameters, indicating that the inverse problem is highly ill-posed. 
This observation is further supported by the numerical experiments presented later.

The inverse problem \eqref{eq:IP-Boltz} needs to be solved as an optimization problem, which requires multiple optimization iterations.
Its objective function needs to be evaluated repeatedly for different values of $\bmu$, requiring the solution of the Boltzmann equation at each optimization iteration.
These are expensive tasks in general, rooted from the intrinsic high dimensionality (hence high-dimensional mesh-based discrete spaces at the numerical level) 
of the phase space on which the Boltzmann equation is defined, as well as the multi-fold integral form of the nonlinear collision operator. Moreover, obtaining steady-state solutions through iterative methods often requires a large number of iterations, particularly in the near fluid regime.   In this work, we are able to greatly improve the computational efficiency of the overall inverse solvers  by working with a low-dimensional representation of the relevant functions through a RB space (as in \Cref{sec:IVP:single-level} for the single-level case) as well as an efficient ROM Boltzmann solver (as in \Cref{sec:IVP:bilevel} for the bilevel case). Both are constructed using a greedy approach and exploit the quadratic structure of the collision operator together with an EIM separable approximation of the collision kernel. 

\subsection{Inverse problem: bilevel optimization}
\label{sec:IVP:bilevel}
Suppose we have obtained the RB space ${\mathbb V}_N=\textrm{span}\{\fFOM(\cdot\,;\bmu_i)\}_{i=1}^N$ and its orthonormal basis  $ \{u_i\}_{i=1}^N$ from \Cref{sec:offline},  we can then use the computationally efficient ROM solver $\text{ROM}(\mathV_N; \bmu)$ to replace the PDE constraint in \eqref{eq:IP-Boltz}.
With this, the inverse problem becomes a bilevel optimization problem, with the lower-level as the constrained optimization \eqref{eq:c-constrained-opt} to minimize the residual to find the reduced solution $\ffROM(\cdot\,; \bmu)=\sum_{i=1}^N c_i(\bmu) \bu_i$, and an implicit objective in the upper level evaluating the data misfit by using this reduced solution:
\begin{equation}\label{eq:IP-Boltz-ROM-bilevel}
\begin{aligned}
\bmu^* = \; & \underset{\bmu=(a,b)}{\text{arg\,min}}\; \;  \left\| T\left(\sum_{i=1}^N c_i(\bmu) \bu_i \right) - T^{\text{data}} \right\|^2_{x}=:J(\bmu)\\
&\text{subject~to } \bc(\bmu)=(c_i(\bmu)) \; \text{solving the constrained optimization \eqref{eq:c-constrained-opt}}.
\end{aligned}
\end{equation}

The bilevel/implicit optimization problem \eqref{eq:IP-Boltz-ROM-bilevel} can be solved by treating the inner-level optimization problem as a black box. More precisely, for any current value of $\bmu$, we solve the constrained optimization problem \eqref{eq:c-constrained-opt} to obtain $\bc(\bmu)$, and then substitute $\bc(\bmu)$ into the outer objective function.
In this way, \eqref{eq:IP-Boltz-ROM-bilevel} can be viewed as an unconstrained optimization problem with an implicit black-box objective evaluator. Therefore, any off-the-shelf unconstrained optimization solver can be applied.  
In our implementation, we employ two approaches to solve this bilevel optimization problem. The first is a gradient-free pattern-search method with a Nelder–Mead search step \cite{audet2002analysis,kolda2003optimization,nelder1965simplex}, implemented using \texttt{patternsearch} with \texttt{SearchFcn} set to \texttt{@searchneldermead} in MATLAB. 
The second is a gradient-based BFGS quasi-Newton method with a cubic line-search procedure \cite{broyden1970convergence,fletcher1970new,goldfarb1970family,shanno1970conditioning}, implemented using \texttt{fminunc} with the \texttt{quasi-newton} algorithm.
The algorithm for solving this inverse problem as a bilevel optimization using either optimization solver is summarized in \Cref{alg:bilevel}.

\begin{remark}
    Note that our ROMs in general do not ensure any smoothness of the reduced coefficient vector $\bc(\bmu)$. This contributes to some non-smoothness of this bilevel optimization problem and hence brings some computational challenge to gradient-based optimization solvers. 
    It is practically relevant for the bilevel optimization strategy to have access of ROMs that generate $\bc(\bmu)$ with certain smoothness, and this will be left to our future investigation. 
\end{remark}

\begin{algorithm}[h]
\caption{Bilevel optimization for inverse problem of Boltzmann equation, using the proposed ROM 
\label{alg:bilevel}}
\begin{algorithmic}[1]

\STATE \textbf{Input:} observation data $\Tdata $, initial parameter $\bmu^{(0)}$, RB basis $\{u_i\}_{i=1}^N$, precomputed terms $\{ \mL (u_i):\ i\in [N]\}$ and  $\{\mathcal{Q}(\bu_i,\bu_j; \mueim_m):\   i,j\in [N],\ m\in[S] \}$, EIM points $\{\mueim_m,r_m\}_{m=1}^S$ and inverse interpolation matrix $\bG^{-1}$

\FOR{$k=0,1,\ldots$}
\STATE \textbf{Lower level:} for $\bmu^{(k)}=(a^{(k)},b^{(k)})$, follow \textbf{EIM Evaluation} \& \textbf{Optimization} steps in \Cref{alg:online} 
to compute the ROM solution $\ffROM(\cdot\,;\bmu^{(k)})$
\STATE \textbf{Upper level:} Evaluate the data misfit objective $ J(\bmu^{(k)}) = \| T(\ffROM(\cdot\,;\bmu^{(k)}))  - T^{\text{data}} \|^2_{x}$; Update $\bmu^{(k+1)} = \bmu^{(k)} + \alpha^{(k)} \Delta \bmu^{(k)}$, where the step size $\alpha^{(k)} $ and the search direction $\Delta \bmu^{(k)}$ are  determined by the current objective value $J(\bmu^{(k)})$ and the chosen optimization algorithm from the off-the-shelf solver for minimizing $J(\bmu)$; terminate if the optimization tolerance is reached 
\ENDFOR
\STATE \textbf{Output:} reconstructed parameter $\bmu^{(k+1)}$
\end{algorithmic}
\end{algorithm}

\subsection{Inverse problem: single-level optimization}
\label{sec:IVP:single-level}

With the low-dimensional ROM solver providing solutions to the Boltzmann equation,  our bilevel optimization above is much more efficient compared with repeatedly calling the expensive FOM solver. On the other hand, the inner level still requires an iterative optimization solve each time $\bmu$ is updated.  

This motivates us to seek an alternative strategy, %
by  reformulating  the task as a single-level optimization problem, in which the lower-level constrained optimization problem is replaced by its Karush--Kuhn--Tucker (KKT) conditions~\cite{guerra2025learning}. The first-order optimality condition for the lower-level optimization \eqref{eq:c-constrained-opt} with the Lagrange multiplier $\lambda\in\mathbb{R}$ is 
$$
\lambda \partial_{c_i} \left[\left\langle \sum_{j=1}^N c_j \bu_j \right\rangle_{x,\bv} - \mass \right]= \partial_{c_i} 
\left\| \Res\left(\sum_{j=1}^N c_j \bu_j; \bmu\right) \right\|_{x, \bv}^2, \, \, i=1,\ldots, N,
$$
where the term on the left is $\lambda \left\langle  \bu_i \right\rangle_{x,\bv}$, and that on the right is nonlinear (indeed, a third-order polynomial) 
in $\bc$ and given in \eqref{eq:dc_Res}.
 Together with the mass conservation constraint (as a linear equation) in \eqref{eq:c-constrained-opt}, these formulate the KKT conditions for \eqref{eq:c-constrained-opt}. Once we use these conditions 
 to replace the lower-level optimization in \eqref{eq:IP-Boltz-ROM-bilevel}, the optimization becomes single-level with nonlinear constraints:
\begin{equation}\label{eq:IP-Boltz-single-level}
\begin{aligned}
\underset{\bmu=(a,b), \bc, \lambda}{\text{minimize}}\; \;  &\left\| T\left(\sum_{i=1}^N c_i \bu_i \right) - T^{\text{data}} \right\|^2_{x}\\
\text{subject~to } &\lambda \left\langle  \bu_i \right\rangle_{x,\bv}=
2 \left\langle
\left[ \sum_{j=1}^N c_j \mL(\bu_j)  - b 
\sum_{k,l=1}^N c_k c_l   \left(\sum_{m=1}^S \beta_m(a) \mathcal{Q}(\bu_k,\bu_l; \mueim_m) \right)
\right] \right. \\
& \;\;
\left.\cdot
\left[ 
\mL( \bu_i)  - b
\sum_{j=1}^N c_j \left( 
\sum_{m=1}^S \beta_m(a) \Big(\mathcal{Q}(\bu_i,\bu_j; \mueim_m) + \mathcal{Q}(\bu_j,\bu_i; \mueim_m) \Big)
\right)
\right] \right\rangle_{x,\bv}, \;   i\in[N] \\
 & \left\langle \sum_{i=1}^N c_i \bu_i \right\rangle_{x,\bv}=\mass.
\end{aligned}
\end{equation}

We want to stress that this single-level optimization problem becomes practically tractable because the low-dimensional RB space $\mathV_N$ is now used, rather than the high-dimensional grid-function space $\mathV$, to represent the relevant functions in the KKT conditions. 
Note that the single-level reformulation \eqref{eq:IP-Boltz-single-level} is not equivalent to the bilevel optimization problem \eqref{eq:IP-Boltz-ROM-bilevel}, and it is a relaxation of the bilevel problem.

The single-level optimization problem \eqref{eq:IP-Boltz-single-level} can be solved using any off-the-shelf nonlinear optimization solvers.
It has a nonlinear objective, which depends only on $\bc$, and nonlinear constraints, which are polynomial in $\bc$, quadratic in $b$, and exponential in $a$. 
Thus, the gradients of both the objective and the constraints can be derived analytically and implemented directly in the optimization solver.
In our implementation, we solve this nonlinear constrained optimization problem using \texttt{fmincon} with the interior-point algorithm \cite{byrd1999interior,byrd2000trust,waltz2006interior} and an L-BFGS Hessian approximation \cite{liu1989limited}.
We do not consider a gradient-free method for solving \eqref{eq:IP-Boltz-single-level}, since the dimension of the unknown variables increases to $3+N$.  This makes it computationally expensive to explore the entire search space without the gradient information. The algorithm for solving our inverse problem as a single-level optimization is summarized in \Cref{alg:single-level}.

\begin{remark}Another advantage of the single-level reformulation is that it allows future extensions to use reduced representations of relevant functions (via low-dimensional RB spaces) in more general inverse problem settings, e.g., where $b$ is not just a constant but a function of $x$. In this case, the computational complexity does not increase significantly with the dimension of the discretization of $b$, 
since the gradient with respect to $b$ can still be derived analytically.
\end{remark}

\begin{algorithm}[h]
\caption{Single-level optimization for inverse problem for Boltzmann equation, using RB space, precomputed collision terms and EIM kernel approximation\label{alg:single-level}}
\begin{algorithmic}[1]

\STATE \textbf{Input:} observation data $\Tdata $, initial parameter $\bmu^{(0)}$ and Lagrange multiplier $\lambda^{(0)}$, RB basis $\{u_i\}_{i=1}^N$, precomputed terms $\{ \mL (u_i):\ i\in [N]\}$ and  $\{\mathcal{Q}(\bu_i,\bu_j; \mueim_m):\   i,j\in [N],\ m\in[S] \}$, EIM points $\{\mueim_m,r_m\}_{m=1}^S$ and inverse interpolation matrix $\bG^{-1}$, QR matrix $\bR$

\STATE Initialize coeff vector $\bc^{(0)} = \bR[1/N, \ldots, 1/N]^\top$.

\FOR{$k=0,1,\ldots$}

\STATE \textbf{Objective evaluation:} compute the objective value $\left\| T\left(\sum_{i=1}^N c^{(k)}_i \bu_i \right) - T^{\text{data}} \right\|^2_{x}$ and its derivative to $\bc$ at current point $\bc^{(k)}$
\STATE \textbf{Constraint evaluation}:
\STATE For $\bmu^{(k)}=(a^{(k)},b^{(k)})$, evaluate $\{\kappa(r_m; a^{(k)})\}_{m=1}^S$ and compute $\bbeta(a^{(k)}) = (\beta_m(a^{(k)}))$ using \eqref{eq:beta-formula} with $\bG^{-1}$
\STATE Compute the values and the gradients of the constraints 
in \eqref{eq:IP-Boltz-single-level} by using $\bbeta(a^{(k)})$ and precomputed terms $\{ \mL (u_i):\ i\in [N]\}$ and  $\{\mathcal{Q}(\bu_i,\bu_j; \mueim_m):\   i,j\in [N],\ m\in[S] \}$

    \STATE \textbf{Update:} update $(\bmu^{(k+1)}, \bc^{(k+1)}, \lambda^{(k+1)})$, using the interior-point method and the current values and gradients of the objective and constraints at $(\bmu^{(k)}, \bc^{(k)}, \lambda^{(k)})$; terminate if  the optimization tolerance is reached
\ENDFOR
\STATE \textbf{Output:} reconstructed parameter $\bmu^{(k+1)}$
\end{algorithmic}
\end{algorithm}

\section{Numerical experiments}
\label{sec:num}

In this section, we conduct several numerical experiments to %
demonstrate the performance of the proposed ROM method in \Cref{sec:ROM} and its application to solving the inverse problem in \Cref{sec:IVP}. 

Our numerical examples are set on the spatial domain $\Omega=[0,1]$,
with boundary temperatures $T_L=1$ and $T_R=2$. The domain in velocity is truncated to $[-L_v, L_v]^2$, with
$L_v = \frac{4.5(3+\sqrt{2})}{2}$. The FOM %
follows the description in \Cref{sec:FOM}, with $N_x=100$ and $N_v=32$, and the FOM solution is computed iteratively until the discrete $L^2$ norm of the residual  \eqref{eq:res} is smaller than $10^{-6}$. 

As discussed in \Cref{sec:EIM}, we first apply EIM to $r\kappa(r;a)=\kappa(r;a+1)$ to generate its separable approximation as in \eqref{eq:sep:approx} with $S=8$ terms for $a\in [a_-, a_+]=[-1,1]$ and $r\in [r_-, r_+]=[0,9]$. This is based on $4000$ equidistant training points in $a$ and in $r$, respectively, with endpoints included (see e.g., p.73 in \cite{hesthaven2016certified}). The resulting %
parameter points  $\{a_m\}_{m=1}^S$ and interpolation points $\{r_m\}_{m=1}^S$ render
an average $L_\infty$ testing error of $5.0E-4$ for $r\kappa(r;a)-r\kappa^S(r;a)$ over $10$ tests, where each test uses
$4000$ random points in $r$ and $500$ points points in $a$, both sampled uniformly. 
While the EIM separable approximations with larger $S$ can be more accurate, the conditioning of the interpolation matrix $\bG$ grows exponentially in $S$ and subsequently the optimization task from the inverse problem in \Cref{sec:num:inverse} becomes more challenging numerically. With $S=8$,  the 2-norm condition number of $\bG$ is about 
$2.78E+4$.

\subsection{Performance of the proposed ROM}
\label{sec:num-rom}

In the following experiments, we study the performance of our proposed ROM for the parametric Boltzmann equation. 
We consider two parameter sets %
corresponding to different regimes: a more 
\emph{collision-dominated} regime with larger $b$, given by $\mM=[-1,1]\times[20,50]$, and a more 
\emph{transport-dominated} regime with smaller $b$, given by $\mM=[-1,1]\times[1,10]$.
The training samples $\Mtrain$ are generated randomly from a uniform distribution over $\mM$, with $|\Mtrain|=50$. 
We perform the ROM algorithm described in \Cref{sec:ROM} until the spectral tolerance reaches $\tau_{\mathrm{tol}} = 3\times 10^{-6}$. 
The constrained optimization problem in \eqref{eq:c-constrained-opt} for finding the reduced solution is solved using MATLAB’s \texttt{fmincon} with the SQP algorithm, analytical objective gradients enabled through \texttt{SpecifyObjectiveGradient}, and optimality tolerance \footnote{The optimality tolerance here is a tolerance for the first-order optimality measure, and only enforces a necessary condition for an optimizer.} set to $10^{-7}$. 
We also evaluate the ROM on a uniformly random testing parameter set $\Mtest\subset \mM$ with $|\Mtest|=20$.

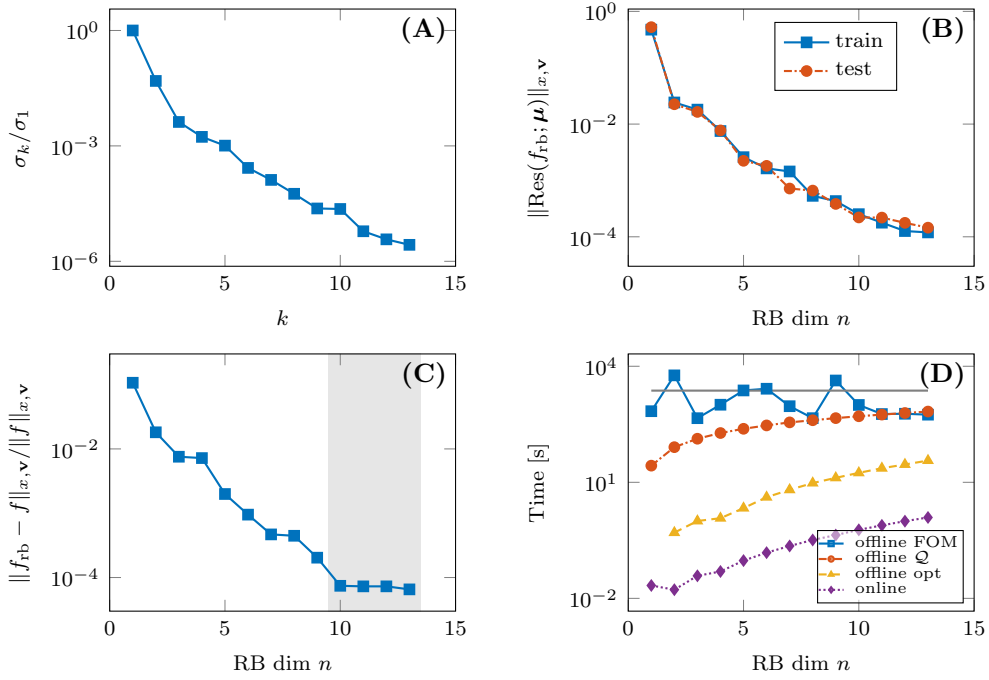
\begin{figure}[tbhp]
    \centering
    \begin{tikzpicture}
    \begin{groupplot}[
        group style={
            group size=2 by 2,
            horizontal sep=65pt,
            vertical sep=33pt,
        },
        width=.42\textwidth,
        height=.34\textwidth,
        cycle list name=siam,
        table/skip first n=2,
        xlabel={RB dim $n$},
        tick label style={font=\scriptsize},
        label style={font=\scriptsize},
        title style={font=\small},
        legend style={font=\scriptsize},
        legend cell align={left},
        xmin=0,
        xmax=15,
        axis on top,
    ]

    \nextgroupplot[
        ylabel={$\sigma_{k}/\sigma_{1}$},
        ymode=log,
        xlabel = {$k$},
    ]
    \addplot table[x=Nmu,y=xiSingularValueNormalized] {data/test-Boltz-qr-fast-eim-S8-UNIa-trunc8-a-11-b2050-restol1e-06-numTrain50-numTest20.txt};

    \nextgroupplot[
        ylabel={$\|\Res(\ffROM;\bmu)\|_{x,\bv}$},
        ymode=log,
    legend style={
        at={(0.6,0.95)},
        anchor=north,
    },
    ymin = 3e-5, ymax =1.1, 
    ]
    \addplot table[x=Nmu,y=trainResidual] {data/test-Boltz-qr-fast-eim-S8-UNIa-trunc8-a-11-b2050-restol1e-06-numTrain50-numTest20.txt};
    \addlegendentry{train}
    \addplot table[x=Nmu,y=maxEIMResidual] {data/test-Boltz-qr-fast-eim-S8-UNIa-trunc8-a-11-b2050-restol1e-06-numTrain50-numTest20.txt};
    \addlegendentry{test}

    \nextgroupplot[
        ylabel={${\|\ffROM - f \|_{x,\bv}}/{\|f \|_{x,\bv}}$},
        ymode=log,
        legend pos=north east,
        xlabel={RB dim $n$},
        ymin = 3e-5, ymax = 3e-1,
    ]
    \addplot[draw=none,name path=shadebottom,forget plot] coordinates {(9.5,1e-10) (13.5,1e-10)};
\addplot[draw=none,name path=shadetop,forget plot] coordinates {(9.5,3) (13.5,3)};
\addplot[fill=gray!20,draw=none,forget plot] fill between[of=shadebottom and shadetop];
    \addplot table[x=Nmu,y=maxRelL2] {data/test-Boltz-qr-fast-eim-S8-UNIa-trunc8-a-11-b2050-restol1e-06-numTrain50-numTest20.txt};

    \nextgroupplot[
        ylabel={Time [s]},
        xlabel={RB dim $n$},
        ymode = log,
        scaled y ticks=false,
        legend columns=1,
        legend pos=south east,
        legend style={
            font=\tiny,
            cells={anchor=west},
            row sep=-4pt,
            column sep=2pt,
            inner xsep=0pt,
            inner ysep=0pt,
        fill=white,
        fill opacity=0.50,
        text opacity=1,
        },
        legend image post style={scale=0.55},
    ]
      \addplot table[x=Nmu,y=traintFOM] {data/test-Boltz-qr-fast-eim-S8-UNIa-trunc8-a-11-b2050-restol1e-06-numTrain50-numTest20.txt};
\addlegendentry{offline FOM}
    \addplot table[x=Nmu,y=traintQ] {data/test-Boltz-qr-fast-eim-S8-UNIa-trunc8-a-11-b2050-restol1e-06-numTrain50-numTest20.txt};
    \addlegendentry{offline $\mathcal{Q}$}
    \addplot table[x=Nmu,y=traintgreedy] {data/test-Boltz-qr-fast-eim-S8-UNIa-trunc8-a-11-b2050-restol1e-06-numTrain50-numTest20.txt};
    \addlegendentry{offline opt}
    \addplot table[x=Nmu,y=meanTestEIMTime] {data/test-Boltz-qr-fast-eim-S8-UNIa-trunc8-a-11-b2050-restol1e-06-numTrain50-numTest20.txt};
    \addlegendentry{online}
    \addplot[black!50, thick] table[x=Nmu,y=meanFOMSolutionTime] {data/test-Boltz-qr-fast-eim-S8-UNIa-trunc8-a-11-b2050-restol1e-06-numTrain50-numTest20.txt};

    \end{groupplot}
    \node[anchor=north east, font=\bfseries\small, inner sep=1pt] at ($(group c1r1.north east)+(-2pt,-2pt)$) {(A)};
    \node[anchor=north east, font=\bfseries\small, inner sep=1pt] at ($(group c2r1.north east)+(-2pt,-2pt)$) {(B)};
    \node[anchor=north east, font=\bfseries\small, inner sep=1pt] at ($(group c1r2.north east)+(-2pt,-2pt)$) {(C)};
    \node[anchor=north east, font=\bfseries\small, inner sep=1pt] at ($(group c2r2.north east)+(-2pt,-2pt)$) {(D)};
    \end{tikzpicture}
    \caption{Training and testing results for the proposed ROM in \emph{collision-dominated} regime: $\bmu \in[-1,1]\times [20,50]$: (A) the normalized singular values of the offline $\rbF_{N}$;  (B) the maximal  residual norms for different RB dimensions, the training residual is maximized over the remaining set $\Mtrain\setminus\{\bmu_i\}_{i=1}^{n}$;
    (C) the maximal relative $L^2$ differences between the FOM and ROM solutions over $\Mtest$, with the plateau effect highlighted in the gray-shaded region;
    (D) the computational time for the (offline) training and (online) testing stages. For the offline stage, we separately record the additional time required for each task when adding one basis vector: computing the FOM solution $f(\cdot\,;\bmu_n)$, precomputing $\mathcal{Q}$ terms, and solving the constrained optimization \eqref{eq:c-constrained-opt} for all the remaining training samples. For the testing stage, we report the average online computation time. The gray
    horizontal line is the average time to compute one FOM solution 
over the test set:
    2345 [s]. } 
    \label{fig:ROM-qr-fast-eim-data}
\end{figure}

We start with  the Boltzmann equation in a relatively collision-dominated regime with $\mM=[-1,1]\times[20,50]$, and report in \Cref{fig:ROM-qr-fast-eim-data} 
the training and testing performance of the ROM algorithm. Upon the termination of the offline greedy procedure, we have an RB space of dimension $N=13$.
Subfigure (A) shows the normalized singular values of $\rbF_{N}$, the matrix representation of the RB space $\mathV_N$, obtained from the offline greedy search algorithm. The singular values decay rapidly in the first three modes, by about three orders of magnitude, and then continue to decay exponentially at a slightly slower rate, decreasing by another three orders of magnitude over the remaining modes.
Subfigure (B) compares the maximum residual norms $\|\Res(\ffROM;\bmu)\|_{x,\bv}$ for the training and testing sets for different dimensions $n$ of the RB space. The training and testing residuals exhibit similar trends: both decay rapidly 
as the RB space grows to two dimensions 
and then continue to decrease at a slightly slower exponential rate. 
The training and testing errors being comparable also evidences that the training parameter set is chosen adequately large.
Subfigure (C) plots the maximum relative discrete $L^2$ error between the FOM and ROM solutions over the test set $\Mtest$, and the overall decay rate stays exponential.
The error plateaus around $n=10$ at approximately $7\times 10^{-5}$ (gray-shaded region)
and does not change much as the RB space is further enriched.
Such a phenomenon can also be observed for ROMs designed for linear kinetic models  \cite{matsuda2025reduced}, and is likely related to the conditioning of the underlying reduced problem.
Based on this observation and (partially) on the offline cost to be discussed next, choosing $n=10$ appears sufficient to represent/approximate the parameter-induced solution manifold.  

In Subfigure (D) of \Cref{fig:ROM-qr-fast-eim-data}, we further record and compare the computational time for the offline training and online testing stages.
\begin{itemize}
    \item 
For the \emph{offline training} stage, we record the additional time required to add one extra dimension to RB space in the ROM greedy algorithm. To get more insights, we separately report the computational time for different components of the algorithm, including solving the $(n-1)$-dimensional constrained optimization problem \eqref{eq:c-constrained-opt} for all remaining $|\Mtrain|-n+1$ training samples, computing one FOM solution $f(\cdot\,;\bmu_n)$, and precomputing the additional $(2n-1)S$ collision terms. It is observed that 
the optimization step used in the greedy search for the next parameter $\bmu_n$ takes relatively much less time than the other components, on the order of tens of seconds. This is because that, after the collision-term precomputation and the EIM kernel approximation, the constrained optimization problem \eqref{eq:c-constrained-opt} becomes an $n$-dimensional polynomial optimization problem with linear constraints, and no longer involves expensive collision-term evaluations.
The offline computational cost is dominated by the FOM solve, which takes several thousand seconds and is at least two orders of magnitude more expensive than the multiple optimization solves.  The fact highlights  the overall efficiency of the MOR framework, especially as one notes that, with the greedy procedure, only $n$ FOM solves are needed in order to build an RB space of dimension $n$.

As the reduced-basis dimension increases, the precomputation cost of the extra $\mathcal{Q}$-related  terms also grows
linearly in $n$ 
and becomes comparable to, or even larger than, the FOM cost when $n>10$. Although this offline cost remains expensive, the EIM kernel approximation \eqref{eq:Q-eim} eliminates the need to compute these collision terms in the online stage. Moreover, even in the offline stage, the precomputation and EIM approximation provide substantial savings: without EIM, one would need $(2n-1)(|\Mtrain|-n)$ expensive $\mathcal{Q}$ evaluations; without precomputation, the collision terms would need to be recomputed in every optimization iteration to find the reduced solution.

\item For the \emph{online testing} stage, we report the average online computation time over the test parameter samples. The gray solid horizontal line indicates the average computation time of the FOM,
which is about $2345$ seconds. In contrast, with the EIM kernel approximation and the $\mathcal{Q}$-term  precomputation, the
ROM online stage only requires solving a single $n$-dimensional polynomial optimization problem.
For example, when $n=10$, one online ROM solve takes less than one second, compared with several thousand seconds for the FOM solve. This yields a substantial speedup of more than $10^3$, while maintaining a relative error of approximately $0.01\%$.
\end{itemize}

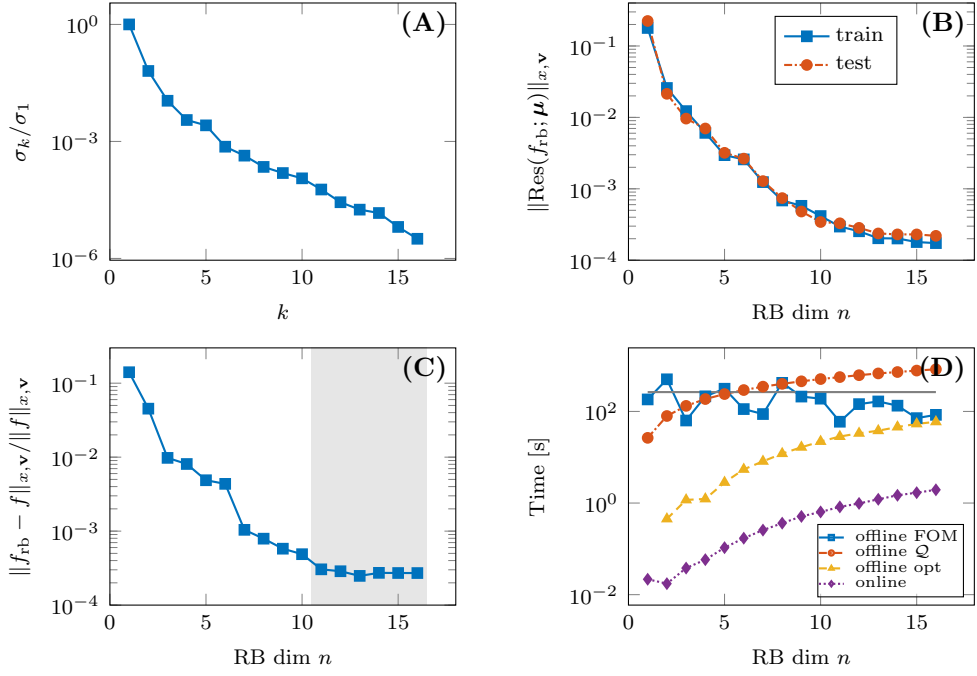
\begin{figure}[h]
    \centering
    \begin{tikzpicture}
    \begin{groupplot}[
        group style={
            group size=2 by 2,
            horizontal sep=65pt,
            vertical sep=33pt,
        },
        width=.42\textwidth,
        height=.34\textwidth,
        cycle list name=siam,
        table/skip first n=2,
        xlabel={RB dim $n$},
        tick label style={font=\scriptsize},
        label style={font=\scriptsize},
        title style={font=\small},
        legend style={font=\scriptsize},
        legend cell align={left},
        xmin=0,
        xmax=18,
        axis on top,
    ]

    \nextgroupplot[
        xlabel = {$k$},
        ylabel={$\sigma_{k}/\sigma_{1}$},
        ymode=log,
    ]
    \addplot table[x=Nmu,y=xiSingularValueNormalized] {data/test-Boltz-qr-fast-eim-S8-UNIa-trunc8-a-11-b110-restol1e-06-numTrain50-numTest20.txt};

    \nextgroupplot[
        ylabel={$\|\Res(\ffROM;\bmu)\|_{x,\bv}$},
        ymode=log,
    legend style={
        at={(0.6,0.95)},
        anchor=north,
    },
    ymin = 1e-4, ymax = 4e-1,
    ]

    \addplot table[x=Nmu,y=trainResidual] {data/test-Boltz-qr-fast-eim-S8-UNIa-trunc8-a-11-b110-restol1e-06-numTrain50-numTest20.txt};
    \addlegendentry{train}
    \addplot table[x=Nmu,y=maxEIMResidual] {data/test-Boltz-qr-fast-eim-S8-UNIa-trunc8-a-11-b110-restol1e-06-numTrain50-numTest20.txt};
    \addlegendentry{test}

    \nextgroupplot[
        ylabel={${\|\ffROM - f \|_{x,\bv}}/{\|f \|_{x,\bv}}$},
        ymode=log,
        legend pos=north east,
        xlabel={RB dim $n$},
        ymin = 1e-4, ymax = 3e-1,
    ]

    \addplot[draw=none,name path=shadebottom,forget plot] coordinates {(10.5,1e-4) (16.5,1e-4)};
\addplot[draw=none,name path=shadetop,forget plot] coordinates {(10.5,3e-1) (16.5,3e-1)};
\addplot[fill=gray!20,draw=none,forget plot] fill between[of=shadebottom and shadetop];

    \addplot table[x=Nmu,y=maxRelL2] {data/test-Boltz-qr-fast-eim-S8-UNIa-trunc8-a-11-b110-restol1e-06-numTrain50-numTest20.txt};

    \nextgroupplot[
        ylabel={Time [s]},
        xlabel={RB dim $n$},
        ymode = log,
        scaled y ticks=false,
        legend columns=1,
        legend pos=south east,
        legend style={
            font=\tiny,
            cells={anchor=west},
            row sep=-4pt,
            column sep=2pt,
            inner xsep=0pt,
            inner ysep=0pt,
            fill=white,
        fill opacity=0.50,
        text opacity=1,
        },
        legend image post style={scale=0.55},
axis on top,
    ]
      \addplot table[x=Nmu,y=traintFOM] {data/test-Boltz-qr-fast-eim-S8-UNIa-trunc8-a-11-b110-restol1e-06-numTrain50-numTest20.txt};
\addlegendentry{offline FOM}
    \addplot table[x=Nmu,y=traintQ] {data/test-Boltz-qr-fast-eim-S8-UNIa-trunc8-a-11-b110-restol1e-06-numTrain50-numTest20.txt};
    \addlegendentry{offline $\mathcal{Q}$}
    \addplot table[x=Nmu,y=traintgreedy] {data/test-Boltz-qr-fast-eim-S8-UNIa-trunc8-a-11-b110-restol1e-06-numTrain50-numTest20.txt};
    \addlegendentry{offline opt}
    \addplot table[x=Nmu,y=meanTestEIMTime] {data/test-Boltz-qr-fast-eim-S8-UNIa-trunc8-a-11-b110-restol1e-06-numTrain50-numTest20.txt};
    \addlegendentry{online}
    \addplot[black!50, thick] table[x=Nmu,y=meanFOMSolutionTime] {data/test-Boltz-qr-fast-eim-S8-UNIa-trunc8-a-11-b110-restol1e-06-numTrain50-numTest20.txt};

    \end{groupplot}
    \node[anchor=north east, font=\bfseries\small, inner sep=1pt] at ($(group c1r1.north east)+(-2pt,-2pt)$) {(A)};
    \node[anchor=north east, font=\bfseries\small, inner sep=1pt] at ($(group c2r1.north east)+(-2pt,-2pt)$) {(B)};
    \node[anchor=north east, font=\bfseries\small, inner sep=1pt] at ($(group c1r2.north east)+(-2pt,-2pt)$) {(C)};
    \node[anchor=north east, font=\bfseries\small, inner sep=1pt] at ($(group c2r2.north east)+(-2pt,-2pt)$) {(D)};
    \end{tikzpicture}
    \caption{Training and testing results for the proposed ROM  \emph{transport-dominated} regime: $\bmu \in[-1,1]\times [1,10]$. The subfigures follow the same descriptions as in \Cref{fig:ROM-qr-fast-eim-data}. The gray horizontal line in (D) is the average time to compute one FOM solution over the test set:
    264 [s].}
    \label{fig:ROM-qr-fast-eim-S8-b110-data}
\end{figure}

We next turn to the experiment  when the Boltzmann equation is in a relatively transport-dominated regime with $\mM=[-1,1]\times[1,10]$, and report the results in \Cref{fig:ROM-qr-fast-eim-S8-b110-data}. Though not rigorously  established, the parameter-induced parameter manifold seems to have a Kolmogorov $n$-width with a relatively slower decay \cite{pinkus2012n,binev2011convergence}, and upon termination of the greedy procedure, a relatively larger RB space of dimension  is obtained. 
Related, the normalized singular values of $\rbF_{N}$ reported in  Subfigure (A)
decay more slowly than those in the more collision-dominated regime shown in \Cref{fig:ROM-qr-fast-eim-data}.
The maximum training and test residual errors in Subfigure (B) and 
the relative $L^2$ errors of the reduced solution over the test set in Subfigure (C)  exhibit decay trends similar to that of the singular values. It is also observed that the relative errors reach a 
plateau from $n\geq 11$, with the plateau value
slightly larger than that in \Cref{fig:ROM-qr-fast-eim-data}, namely approximately $2.7\times 10^{-4}$ compared with $7\times 10^{-5}$.
For the computational time shown in Subfigure (D), the FOM requires fewer iterations to converge and therefore has a lower average computational time of approximately $264$ seconds. In the offline stage, the dominant costs come from the FOM solves and the precomputation of the $\mathcal{Q}$ terms, with the latter becoming more dominant when $n>5$. Finally, the online ROM computation  remains substantially faster than the FOM, achieving a speedup by a factor of at least $200$.

\subsection{Inverse problem} %
\label{sec:num:inverse}

In this section, we apply the ROM trained in \Cref{sec:num-rom} to the inverse problem \eqref{eq:IP-Boltz}. The same setup in \Cref{sec:num-rom} is used for the FOM solution, the EIM approximation and the constrained optimization solver to compute the ROM solution in \eqref{eq:c-constrained-opt}. We use the following three optimization approaches  discussed in \Cref{sec:IVP} to solve the inverse problem:
(1) the bilevel optimization method in \Cref{alg:bilevel} with the gradient-free solver \texttt{patternsearch};
(2) the bilevel optimization method in \Cref{alg:bilevel} with the gradient-based quasi-Newton solver \texttt{fminunc}; and
(3) the single-level optimization formulation \eqref{eq:IP-Boltz-single-level} with the gradient-based interior-point method implemented in \texttt{fmincon}.
We study and test these approaches for the case with noiseless data, where the observation data for the temperature $T$ are generated exactly from the FOM solution corresponding to the true parameter $\bmu^*=(a^*, b^*)$ in both the collision- and transport- dominated regimes.  
Comparison is  further performed between the proposed ROM-based inverse solvers with that based on the FOM. 

Based on our numerical experiments, the following settings are taken  for these optimization approaches.
For the gradient-free bilevel approach, we use \texttt{patternsearch} with mesh tolerance $10^{-4}$, maximum iteration number $1000$, mesh expansion and contraction factors $2$ and $0.5$, respectively, mesh scaling enabled, and the Nelder–Mead search function \texttt{@searchneldermead}. For the gradient-based bilevel approach, we use \texttt{fminunc} with the quasi-Newton algorithm, optimality tolerance $10^{-8}$, step tolerance $10^{-11}$, and maximum number of function evaluations $10^3$. For the single-level formulation, we use \texttt{fmincon} with the interior-point algorithm and the L-BFGS Hessian approximation. The optimality tolerance and step tolerance are set to $10^{-8}$ and $10^{-12}$, respectively, the maximum number of function evaluations is set to $10^6$, the problem is scaled with respect to both the objective and constraints, and analytic gradients are supplied for both the objective and constraints.
We directly impose the lower and upper bounds of $\mM$ in \texttt{fmincon}. For \texttt{patternsearch} and \texttt{fminunc}, we enforce the parameter range by penalization: if a parameter value lies outside $\mM$, the objective function is assigned a large value, set to $10^{10}$.

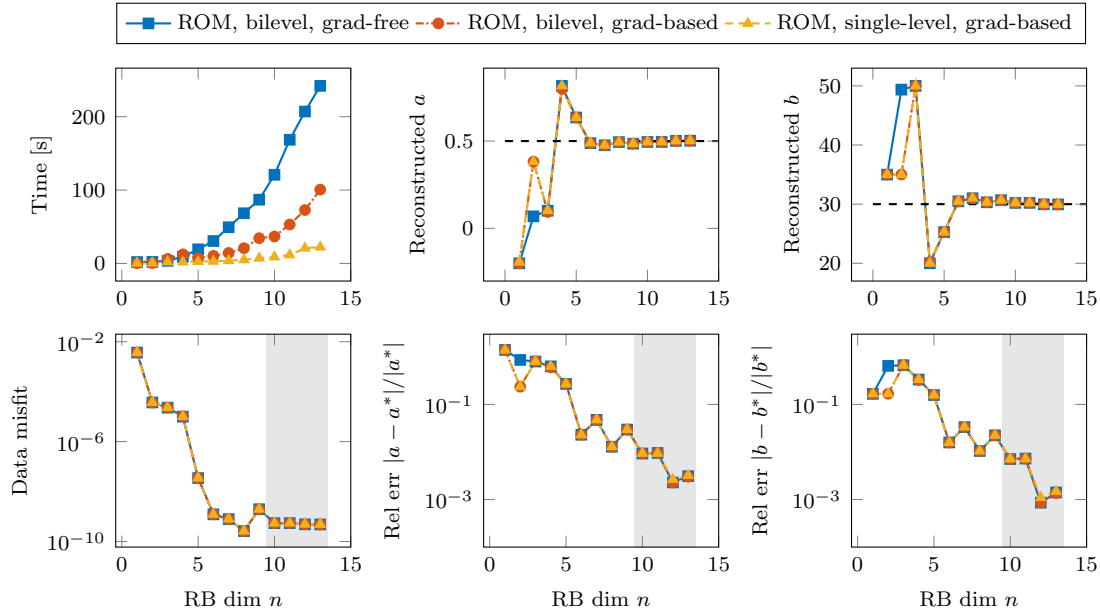
\begin{figure}[h]
	\centering
	\begin{tikzpicture}
	\begin{groupplot}[
	group style={
		group size=3 by 2,
		horizontal sep=50pt,
		vertical sep =20pt,
	},
	width=.32\textwidth,
	height=.3\textwidth,
	legend style={font=\scriptsize},
	tick label style={font=\scriptsize},
	label style={font=\scriptsize},
	xlabel = {RB dim $n$},
	legend style={font=\scriptsize},
	tick label style={font=\scriptsize},
	xmax = 15,
    cycle list name=siam,
    axis on top,
	]
	
	\nextgroupplot[
	xlabel={},
	ylabel={Time [s]},
	legend style = {nodes=right, legend to name=grouplegend-ivp-Boltz-positive-a-S8},
	legend columns=4, 
	]
	\addplot table[x=k,y = t] {data/opt-T0Boltz-qr-fast-eim-S8-UNIa-trunc8-a-11-b2050-restol1e-06-numTrain50-atrue0.5-btrue30-noise0.txt};
	\addlegendentry{ROM, bilevel, grad-free}
	\addplot table[x=k,y = t] {data/opt-T1Boltz-qr-fast-eim-S8-UNIa-trunc8-a-11-b2050-restol1e-06-numTrain50-atrue0.5-btrue30-noise0.txt};
	\addlegendentry{ROM, bilevel, grad-based}
   	\addplot table[x=k,y = t] {data/optKKT-TBoltz-qr-fast-eim-S8-UNIa-trunc8-a-11-b2050-restol1e-06-numTrain50-atrue0.5-btrue30-noise0.txt};
	\addlegendentry{ROM, single-level, grad-based}

	\nextgroupplot[
	ylabel={Reconstructed $a$},
	xlabel={},
	]

    \addplot table[x=k,y = a] {data/opt-T0Boltz-qr-fast-eim-S8-UNIa-trunc8-a-11-b2050-restol1e-06-numTrain50-atrue0.5-btrue30-noise0.txt};
	\addplot table[x=k,y = a] {data/opt-T1Boltz-qr-fast-eim-S8-UNIa-trunc8-a-11-b2050-restol1e-06-numTrain50-atrue0.5-btrue30-noise0.txt};
   	\addplot table[x=k,y = a] {data/optKKT-TBoltz-qr-fast-eim-S8-UNIa-trunc8-a-11-b2050-restol1e-06-numTrain50-atrue0.5-btrue30-noise0.txt};
    \addplot[black, dashed, thick,domain=0:15,] {0.5};

	\nextgroupplot[
	ylabel={Reconstructed $b$},
	xlabel={},
	]
     \addplot table[x=k,y = b] {data/opt-T0Boltz-qr-fast-eim-S8-UNIa-trunc8-a-11-b2050-restol1e-06-numTrain50-atrue0.5-btrue30-noise0.txt};
	\addplot table[x=k,y = b] {data/opt-T1Boltz-qr-fast-eim-S8-UNIa-trunc8-a-11-b2050-restol1e-06-numTrain50-atrue0.5-btrue30-noise0.txt};
   	\addplot table[x=k,y = b] {data/optKKT-TBoltz-qr-fast-eim-S8-UNIa-trunc8-a-11-b2050-restol1e-06-numTrain50-atrue0.5-btrue30-noise0.txt};
    \addplot[black, dashed, thick,domain=0:15,] {30};

	\nextgroupplot[
	ylabel={Data misfit},
	ymode=log,
    ymin = 6e-11, ymax = 2e-2,
	]
        \addplot[draw=none,name path=shadebottom,forget plot] coordinates {(9.5,1e-11) (13.5,1e-11)};
\addplot[draw=none,name path=shadetop,forget plot] coordinates {(9.5,3) (13.5,3)};
\addplot[fill=gray!20,draw=none,forget plot] fill between[of=shadebottom and shadetop];
	     \addplot table[x=k,y = fval] {data/opt-T0Boltz-qr-fast-eim-S8-UNIa-trunc8-a-11-b2050-restol1e-06-numTrain50-atrue0.5-btrue30-noise0.txt};
	\addplot table[x=k,y = fval] {data/opt-T1Boltz-qr-fast-eim-S8-UNIa-trunc8-a-11-b2050-restol1e-06-numTrain50-atrue0.5-btrue30-noise0.txt};
   	\addplot table[x=k,y = fval] {data/optKKT-TBoltz-qr-fast-eim-S8-UNIa-trunc8-a-11-b2050-restol1e-06-numTrain50-atrue0.5-btrue30-noise0.txt};

	\nextgroupplot[
	ylabel={Rel err $|a-a^*|/|a^*|$},
	ymode=log,
      ymin = 1e-4, ymax=3,
	]
            \addplot[draw=none,name path=shadebottom,forget plot] coordinates {(9.5,1e-11) (13.5,1e-11)};
\addplot[draw=none,name path=shadetop,forget plot] coordinates {(9.5,3) (13.5,3)};
\addplot[fill=gray!20,draw=none,forget plot] fill between[of=shadebottom and shadetop];
       \addplot table[x=k,y = relEa] {data/opt-T0Boltz-qr-fast-eim-S8-UNIa-trunc8-a-11-b2050-restol1e-06-numTrain50-atrue0.5-btrue30-noise0.txt};
	\addplot table[x=k,y = relEa] {data/opt-T1Boltz-qr-fast-eim-S8-UNIa-trunc8-a-11-b2050-restol1e-06-numTrain50-atrue0.5-btrue30-noise0.txt};
   	\addplot table[x=k,y = relEa] {data/optKKT-TBoltz-qr-fast-eim-S8-UNIa-trunc8-a-11-b2050-restol1e-06-numTrain50-atrue0.5-btrue30-noise0.txt};

	\nextgroupplot[
	ylabel={Rel err $|b-b^*|/|b^*|$},
	ymode=log,
       ymin = 1e-4, ymax=3,
	]
    \addplot[draw=none,name path=shadebottom,forget plot] coordinates {(9.5,1e-11) (13.5,1e-11)};
\addplot[draw=none,name path=shadetop,forget plot] coordinates {(9.5,3) (13.5,3)};
\addplot[fill=gray!20,draw=none,forget plot] fill between[of=shadebottom and shadetop];
        \addplot table[x=k,y = relEb] {data/opt-T0Boltz-qr-fast-eim-S8-UNIa-trunc8-a-11-b2050-restol1e-06-numTrain50-atrue0.5-btrue30-noise0.txt};
	\addplot table[x=k,y = relEb] {data/opt-T1Boltz-qr-fast-eim-S8-UNIa-trunc8-a-11-b2050-restol1e-06-numTrain50-atrue0.5-btrue30-noise0.txt};
   	\addplot table[x=k,y = relEb] {data/optKKT-TBoltz-qr-fast-eim-S8-UNIa-trunc8-a-11-b2050-restol1e-06-numTrain50-atrue0.5-btrue30-noise0.txt};
	
	\end{groupplot}
	\node[black] at ($(group c2r1) + (0cm,2cm)$) {\pgfplotslegendfromname{grouplegend-ivp-Boltz-positive-a-S8}}; 
	\end{tikzpicture}
	\caption{Inverse problem in \emph{collision-dominated} regime: Comparison of ROM-based inverse problem solvers using 
    different optimization approaches. The true parameters are $a^*=0.5$ and $b^*=30$, indicated by the black dashed lines, and the initial optimization point is set to $a=-0.2$ and $b=35$. 
    The ROM uses the RB spaces (of different dimensions) trained in \Cref{fig:ROM-qr-fast-eim-data} of \Cref{sec:num-rom}.
    We report the total computational time, data misfit (the final objective value in \eqref{eq:IP-Boltz}), the reconstructed values of $a$ and $b$, and their relative errors. 
   The gray-shaded regions in subfigures refer to the RB dimensions for which a plateau is observed in the testing results (C) in \Cref{fig:ROM-qr-fast-eim-data}.
	}\label{fig:ivp-Boltz-positive-a}
\end{figure}

\begin{figure}[h]
	\centering
	\begin{tikzpicture}
	\begin{groupplot}[
	group style={
		group size=3 by 2,
		horizontal sep=50pt,
		vertical sep =20pt,
	},
	width=.32\textwidth,
	height=.3\textwidth,
	legend style={font=\scriptsize},
	tick label style={font=\scriptsize},
	label style={font=\scriptsize},
	xlabel = {RB dim $n$},
	legend style={font=\scriptsize},
	tick label style={font=\scriptsize},
	xmax = 15,
    cycle list name=siam,
    axis on top,
	]
	
	\nextgroupplot[
	xlabel={},
	ylabel={Time [s]},
	legend style = {nodes=right, legend to name=grouplegend-ivp-Boltz-negative-a-S8},
	legend columns=4, 
	]
	\addplot table[x=k,y = t] {data/opt-T0Boltz-qr-fast-eim-S8-UNIa-trunc8-a-11-b2050-restol1e-06-numTrain50-atrue-0.5-btrue30-noise0.txt};
	\addlegendentry{ROM, bilevel, grad-free}
	\addplot table[x=k,y = t] {data/opt-T1Boltz-qr-fast-eim-S8-UNIa-trunc8-a-11-b2050-restol1e-06-numTrain50-atrue-0.5-btrue30-noise0.txt};
	\addlegendentry{ROM, bilevel, grad-based}
   	\addplot table[x=k,y = t] {data/optKKT-TBoltz-qr-fast-eim-S8-UNIa-trunc8-a-11-b2050-restol1e-06-numTrain50-atrue-0.5-btrue30-noise0.txt};
	\addlegendentry{ROM, single-level, grad-based}

	\nextgroupplot[
	ylabel={Reconstructed $a$},
	xlabel={},
	]

    \addplot table[x=k,y = a] {data/opt-T0Boltz-qr-fast-eim-S8-UNIa-trunc8-a-11-b2050-restol1e-06-numTrain50-atrue-0.5-btrue30-noise0.txt};
	\addplot table[x=k,y = a] {data/opt-T1Boltz-qr-fast-eim-S8-UNIa-trunc8-a-11-b2050-restol1e-06-numTrain50-atrue-0.5-btrue30-noise0.txt};
   	\addplot table[x=k,y = a] {data/optKKT-TBoltz-qr-fast-eim-S8-UNIa-trunc8-a-11-b2050-restol1e-06-numTrain50-atrue-0.5-btrue30-noise0.txt};
    \addplot[black, dashed, thick,domain=0:15,] {-0.5};

	\nextgroupplot[
	ylabel={Reconstructed $b$},
	xlabel={},
	]
     \addplot table[x=k,y = b] {data/opt-T0Boltz-qr-fast-eim-S8-UNIa-trunc8-a-11-b2050-restol1e-06-numTrain50-atrue-0.5-btrue30-noise0.txt};
	\addplot table[x=k,y = b] {data/opt-T1Boltz-qr-fast-eim-S8-UNIa-trunc8-a-11-b2050-restol1e-06-numTrain50-atrue-0.5-btrue30-noise0.txt};
   	\addplot table[x=k,y = b] {data/optKKT-TBoltz-qr-fast-eim-S8-UNIa-trunc8-a-11-b2050-restol1e-06-numTrain50-atrue-0.5-btrue30-noise0.txt};
    \addplot[black, dashed, thick,domain=0:15,] {30};

	\nextgroupplot[
	ylabel={Data misfit},
	ymode=log,
    ymin = 4e-11, ymax = 1.5e-4,
	]
        \addplot[draw=none,name path=shadebottom,forget plot] coordinates {(9.5,1e-11) (13.5,1e-11)};
\addplot[draw=none,name path=shadetop,forget plot] coordinates {(9.5,3) (13.5,3)};
\addplot[fill=gray!20,draw=none,forget plot] fill between[of=shadebottom and shadetop];
	     \addplot table[x=k,y = fval] {data/opt-T0Boltz-qr-fast-eim-S8-UNIa-trunc8-a-11-b2050-restol1e-06-numTrain50-atrue-0.5-btrue30-noise0.txt};
	\addplot table[x=k,y = fval] {data/opt-T1Boltz-qr-fast-eim-S8-UNIa-trunc8-a-11-b2050-restol1e-06-numTrain50-atrue-0.5-btrue30-noise0.txt};
   	\addplot table[x=k,y = fval] {data/optKKT-TBoltz-qr-fast-eim-S8-UNIa-trunc8-a-11-b2050-restol1e-06-numTrain50-atrue-0.5-btrue30-noise0.txt};

	\nextgroupplot[
	ylabel={Rel err $|a-a^*|/|a^*|$},
	ymode=log,
     ymin = 1e-4, ymax=3,
	]
        \addplot[draw=none,name path=shadebottom,forget plot] coordinates {(9.5,1e-11) (13.5,1e-11)};
\addplot[draw=none,name path=shadetop,forget plot] coordinates {(9.5,3) (13.5,3)};
\addplot[fill=gray!20,draw=none,forget plot] fill between[of=shadebottom and shadetop];
       \addplot table[x=k,y = relEa] {data/opt-T0Boltz-qr-fast-eim-S8-UNIa-trunc8-a-11-b2050-restol1e-06-numTrain50-atrue-0.5-btrue30-noise0.txt};
	\addplot table[x=k,y = relEa] {data/opt-T1Boltz-qr-fast-eim-S8-UNIa-trunc8-a-11-b2050-restol1e-06-numTrain50-atrue-0.5-btrue30-noise0.txt};
   	\addplot table[x=k,y = relEa] {data/optKKT-TBoltz-qr-fast-eim-S8-UNIa-trunc8-a-11-b2050-restol1e-06-numTrain50-atrue-0.5-btrue30-noise0.txt};

	\nextgroupplot[
	ylabel={Rel err $|b-b^*|/|b^*|$},
	ymode=log,
     ymin = 1e-4, ymax=3,
	]
        \addplot[draw=none,name path=shadebottom,forget plot] coordinates {(9.5,1e-11) (13.5,1e-11)};
\addplot[draw=none,name path=shadetop,forget plot] coordinates {(9.5,3) (13.5,3)};
\addplot[fill=gray!20,draw=none,forget plot] fill between[of=shadebottom and shadetop];
        \addplot table[x=k,y = relEb] {data/opt-T0Boltz-qr-fast-eim-S8-UNIa-trunc8-a-11-b2050-restol1e-06-numTrain50-atrue-0.5-btrue30-noise0.txt};
	\addplot table[x=k,y = relEb] {data/opt-T1Boltz-qr-fast-eim-S8-UNIa-trunc8-a-11-b2050-restol1e-06-numTrain50-atrue-0.5-btrue30-noise0.txt};
   	\addplot table[x=k,y = relEb] {data/optKKT-TBoltz-qr-fast-eim-S8-UNIa-trunc8-a-11-b2050-restol1e-06-numTrain50-atrue-0.5-btrue30-noise0.txt};
	
	\end{groupplot}
	\node[black] at ($(group c2r1) + (0cm,2cm)$) {\pgfplotslegendfromname{grouplegend-ivp-Boltz-negative-a-S8}};
	\end{tikzpicture}
	\caption{Inverse problem in \emph{collision-dominated} regime: Comparison of ROM-based inverse problem solvers using 
    different optimization approaches. The true parameters are $a^*=-0.5$ and $b^*=30$, indicated by the black dashed lines, other setups are the same as in \Cref{fig:ivp-Boltz-positive-a}.
	}\label{fig:ivp-Boltz-negative-a-1}
\end{figure}
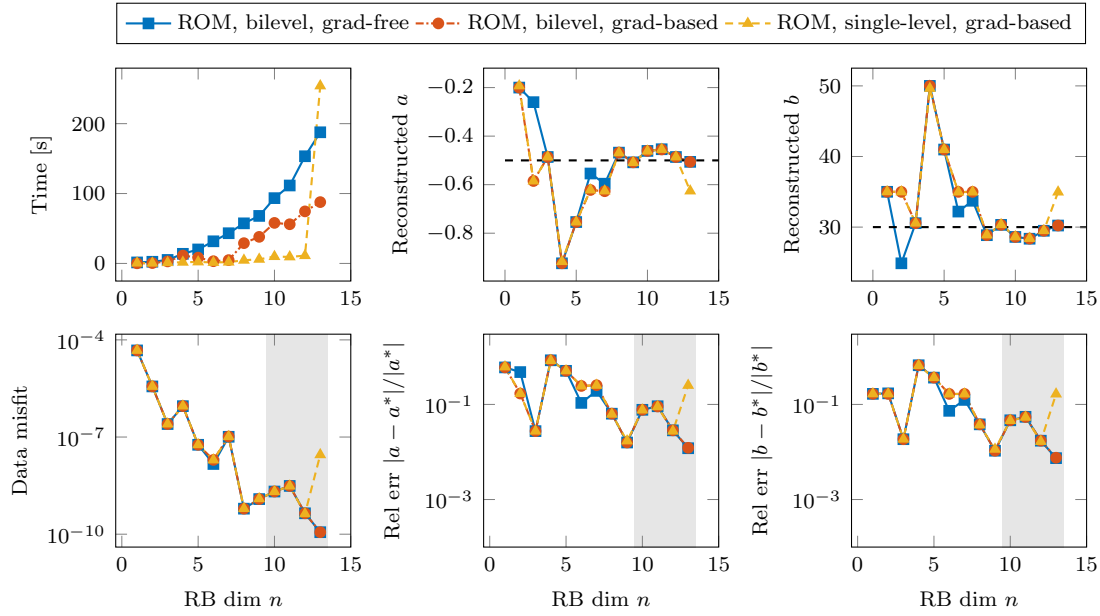
\begin{figure}[h]
	\centering
	\begin{tikzpicture}
	\begin{groupplot}[
	group style={
		group size=3 by 2,
		horizontal sep=50pt,
		vertical sep =20pt,
	},
	width=.32\textwidth,
	height=.3\textwidth,
	legend style={font=\scriptsize},
	tick label style={font=\scriptsize},
	label style={font=\scriptsize},
	xlabel = {RB dim $n$},
	xmax = 15,
    cycle list name=siam,
     axis on top,
	]
	
	\nextgroupplot[
	xlabel={},
	ylabel={Time [s]},
	legend style = {nodes=right, legend to name=grouplegend-ivp-Boltz-negative-a-2-S8},
	legend columns=4,
	]
	\addplot table[x=k,y = t] {data/opt-T0Boltz-qr-fast-eim-S8-UNIa-trunc8-a-11-b2050-restol1e-06-numTrain50-atrue-0.6-btrue25-noise0.txt};
	\addlegendentry{ROM, bilevel, grad-free}
	\addplot table[x=k,y = t] {data/opt-T1Boltz-qr-fast-eim-S8-UNIa-trunc8-a-11-b2050-restol1e-06-numTrain50-atrue-0.6-btrue25-noise0.txt};
	\addlegendentry{ROM, bilevel, grad-based}
    \addplot table[x=k,y = t] {data/optKKT-TBoltz-qr-fast-eim-S8-UNIa-trunc8-a-11-b2050-restol1e-06-numTrain50-atrue-0.6-btrue25-noise0.txt};
	\addlegendentry{ROM, single-level, grad-based}

	\nextgroupplot[
	ylabel={Reconstructed $a$},
	xlabel={},
	]
    \addplot table[x=k,y = a] {data/opt-T0Boltz-qr-fast-eim-S8-UNIa-trunc8-a-11-b2050-restol1e-06-numTrain50-atrue-0.6-btrue25-noise0.txt};
	\addplot table[x=k,y = a] {data/opt-T1Boltz-qr-fast-eim-S8-UNIa-trunc8-a-11-b2050-restol1e-06-numTrain50-atrue-0.6-btrue25-noise0.txt};
    \addplot table[x=k,y = a] {data/optKKT-TBoltz-qr-fast-eim-S8-UNIa-trunc8-a-11-b2050-restol1e-06-numTrain50-atrue-0.6-btrue25-noise0.txt};
    \addplot[black, dashed, thick,domain=0:15,] {-0.6};

	\nextgroupplot[
	ylabel={Reconstructed $b$},
	xlabel={},
	]
    \addplot table[x=k,y = b] {data/opt-T0Boltz-qr-fast-eim-S8-UNIa-trunc8-a-11-b2050-restol1e-06-numTrain50-atrue-0.6-btrue25-noise0.txt};
	\addplot table[x=k,y = b] {data/opt-T1Boltz-qr-fast-eim-S8-UNIa-trunc8-a-11-b2050-restol1e-06-numTrain50-atrue-0.6-btrue25-noise0.txt};
    \addplot table[x=k,y = b] {data/optKKT-TBoltz-qr-fast-eim-S8-UNIa-trunc8-a-11-b2050-restol1e-06-numTrain50-atrue-0.6-btrue25-noise0.txt};
    \addplot[black, dashed, thick,domain=0:15,] {25};

	\nextgroupplot[
	ylabel={Data misfit},
	ymode=log,
    ymin = 3e-11, ymax = 7e-3,
	]
            \addplot[draw=none,name path=shadebottom,forget plot] coordinates {(9.5,1e-11) (13.5,1e-11)};
\addplot[draw=none,name path=shadetop,forget plot] coordinates {(9.5,3) (13.5,3)};
\addplot[fill=gray!20,draw=none,forget plot] fill between[of=shadebottom and shadetop];
    \addplot table[x=k,y = fval] {data/opt-T0Boltz-qr-fast-eim-S8-UNIa-trunc8-a-11-b2050-restol1e-06-numTrain50-atrue-0.6-btrue25-noise0.txt};
	\addplot table[x=k,y = fval] {data/opt-T1Boltz-qr-fast-eim-S8-UNIa-trunc8-a-11-b2050-restol1e-06-numTrain50-atrue-0.6-btrue25-noise0.txt};
    \addplot table[x=k,y = fval] {data/optKKT-TBoltz-qr-fast-eim-S8-UNIa-trunc8-a-11-b2050-restol1e-06-numTrain50-atrue-0.6-btrue25-noise0.txt};

	\nextgroupplot[
	ylabel={Rel err $|a-a^*|/|a^*|$},
	ymode=log,
     ymin = 1e-4, ymax=3,
	]
            \addplot[draw=none,name path=shadebottom,forget plot] coordinates {(9.5,1e-11) (13.5,1e-11)};
\addplot[draw=none,name path=shadetop,forget plot] coordinates {(9.5,3) (13.5,3)};
\addplot[fill=gray!20,draw=none,forget plot] fill between[of=shadebottom and shadetop];
    \addplot table[x=k,y = relEa] {data/opt-T0Boltz-qr-fast-eim-S8-UNIa-trunc8-a-11-b2050-restol1e-06-numTrain50-atrue-0.6-btrue25-noise0.txt};
	\addplot table[x=k,y = relEa] {data/opt-T1Boltz-qr-fast-eim-S8-UNIa-trunc8-a-11-b2050-restol1e-06-numTrain50-atrue-0.6-btrue25-noise0.txt};
    \addplot table[x=k,y = relEa] {data/optKKT-TBoltz-qr-fast-eim-S8-UNIa-trunc8-a-11-b2050-restol1e-06-numTrain50-atrue-0.6-btrue25-noise0.txt};

	\nextgroupplot[
	ylabel={Rel err $|b-b^*|/|b^*|$},
	ymode=log,
     ymin = 1e-4, ymax=3,
	]
            \addplot[draw=none,name path=shadebottom,forget plot] coordinates {(9.5,1e-11) (13.5,1e-11)};
\addplot[draw=none,name path=shadetop,forget plot] coordinates {(9.5,3) (13.5,3)};
\addplot[fill=gray!20,draw=none,forget plot] fill between[of=shadebottom and shadetop];
    \addplot table[x=k,y = relEb] {data/opt-T0Boltz-qr-fast-eim-S8-UNIa-trunc8-a-11-b2050-restol1e-06-numTrain50-atrue-0.6-btrue25-noise0.txt};
	\addplot table[x=k,y = relEb] {data/opt-T1Boltz-qr-fast-eim-S8-UNIa-trunc8-a-11-b2050-restol1e-06-numTrain50-atrue-0.6-btrue25-noise0.txt};
    \addplot table[x=k,y = relEb] {data/optKKT-TBoltz-qr-fast-eim-S8-UNIa-trunc8-a-11-b2050-restol1e-06-numTrain50-atrue-0.6-btrue25-noise0.txt};
	
	\end{groupplot}
	\node[black] at ($(group c2r1) + (0cm,2cm)$) {\pgfplotslegendfromname{grouplegend-ivp-Boltz-negative-a-2-S8}};
	\end{tikzpicture}
	\caption{Inverse problem in \emph{collision-dominated} regime: Comparison of ROM-based inverse problem solvers using 
    different optimization approaches. The true parameters are $a^*=-0.6$ and $b^*=25$, indicated by the black dashed lines, other setups are the same as in \Cref{fig:ivp-Boltz-positive-a}.
	}\label{fig:ivp-Boltz-negative-a-2}
\end{figure}
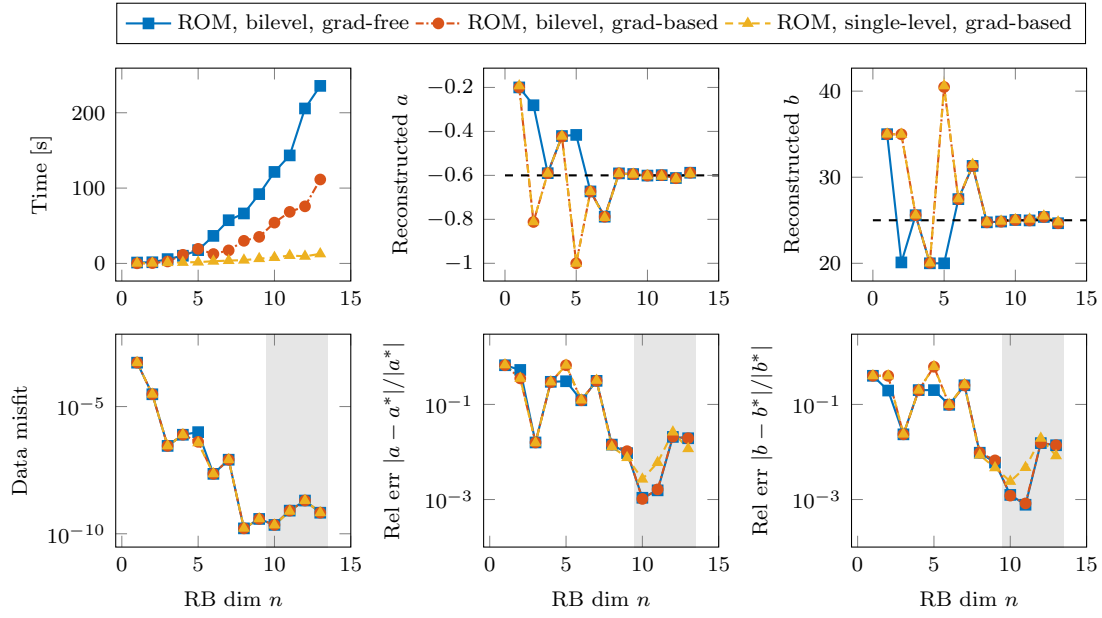

\medskip
\noindent {\bf Inverse problem in collision-dominated regime.}
By using the ROM trained for the relatively collision-dominated regime $\mM=[-1,1]\times[20,50]$,
we first test the three optimization approaches to solve the inverse problem, when the true parameter $\bmu^*=(a^*, b^*)$ is one of the following: $(0.5, 30)$, $(-0.5, 30)$, $(-0.6, 25)$.
For all cases, the initial guess of the parameter is 
set with $\bmu=(-0.2, 35)$. 
The results are shown in \Cref{fig:ivp-Boltz-positive-a,fig:ivp-Boltz-negative-a-1,fig:ivp-Boltz-negative-a-2}. %

Recall from the previous section, the ROM solver we built becomes more accurate as the RB space is enriched and its dimension $n$ grows, until $n$ reaches $10$ and the testing errors of the ROMs display a plateau effect (see \Cref{fig:ROM-qr-fast-eim-data}),  likely due to the conditioning of the ROM.  When the reduced solver is used in the three optimization approaches to solve the inverse problem, the first thing we want to examine is how the inferred parameter is improved as $n$ increases.  %
In all cases, the data misfit shows an overall decreasing trend until $n=10$ in \Cref{fig:ivp-Boltz-positive-a} or $n=8$ in \Cref{fig:ivp-Boltz-negative-a-1,fig:ivp-Boltz-negative-a-2}, after which it either levels off or increases. When we come to the reconstruction errors, 
for the case with the positive $a^* = 0.5$, all three approaches perform as one would hope, in that the relative reconstruction errors for both $a$ and $b$ generally decrease with $n$. They are about $1\%$ when $n = 10$ and about $0.1\%$ when $n = 12$.
For the case with the negative  $a^*=-0.5, -0.6$, the relative reconstruction errors start decreasing roughly from  $n=5$~\footnote{where the ROM reaches a relative testing error of $1.99\times 10^{-3}$.
}, and they then increase from or beyond $n=10$ yet stay below $10\%$. It is a good balance of accuracy and robustness if one uses the ROM near the onset of the plateau (with $n=9$ or $10$) in our inverse problem, with the corresponding reconstruction errors approximately $1\%$ in \Cref{fig:ivp-Boltz-negative-a-1} and $0.1\%$ (or even smaller) in \Cref{fig:ivp-Boltz-negative-a-2}.

As expected, for each optimization approach, the computational time grows with the RB dimension $n$. With a fixed RB space, the three approaches differ quite much in their computational time. For example, when $n=10$, the gradient-free bilevel method requires about $100$ seconds to converge, while the gradient-based bilevel method reduces the computational time by approximately %
half. The single-level formulation provides the most significant acceleration, requiring less than $10$ seconds.
These computational times further validate the use of the ROM to solve the inverse problem.
Recall from the study in \Cref{sec:num-rom},  one FOM solve requires about $2000$ seconds on average. In an inverse problem, where %
forward solves are needed repeatedly during the optimization process, directly calling the FOM would therefore be prohibitively expensive, while ROM brings down the cost to just  several seconds.
The only exception occurs in \Cref{fig:ivp-Boltz-negative-a-1} when $n=13$, where the single-level method takes more time than the gradient-free bilevel method. In this case, the reconstruction errors are also larger and the data misfit is relatively high, suggesting that the optimization solver is likely trapped in a local minimum.

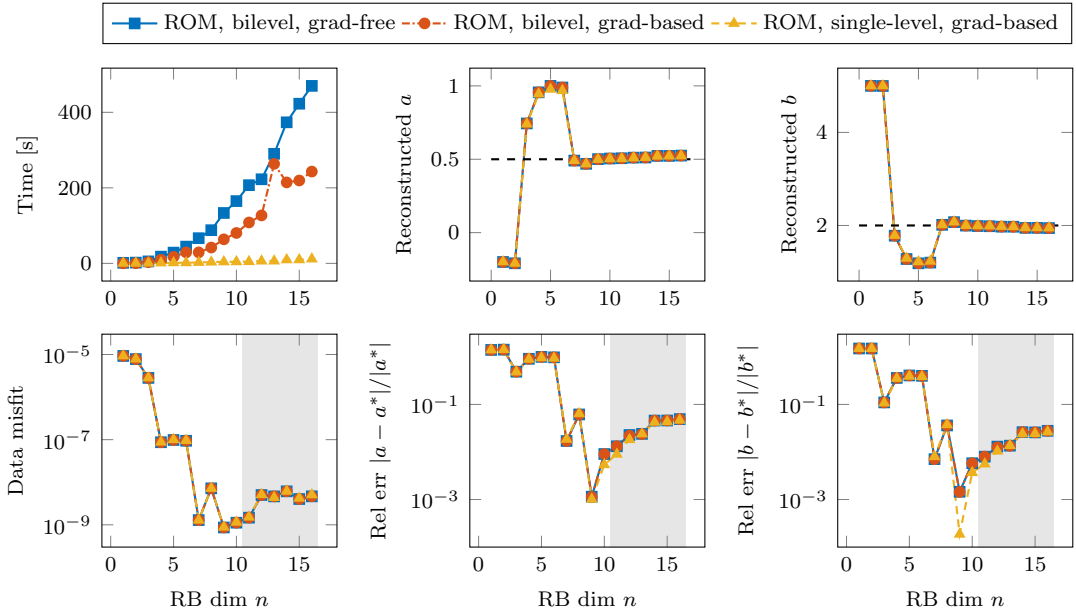
\begin{figure}[h]
	\centering
	\begin{tikzpicture}
	\begin{groupplot}[
	group style={
		group size=3 by 2,
		horizontal sep=50pt,
		vertical sep =20pt,
	},
	width=.32\textwidth,
	height=.3\textwidth,
	legend style={font=\scriptsize},
	tick label style={font=\scriptsize},
		label style={font=\scriptsize},
	xlabel = {RB dim $n$},
	legend style={font=\scriptsize},
	tick label style={font=\scriptsize},
	xmax = 18,
    cycle list name=siam,
    axis on top,
	]
	
	\nextgroupplot[
	xlabel={},
	ylabel={Time [s]},
	legend style = {nodes=right, legend to name=grouplegend-ivp-Boltz-b110-positive-a-S8},
	legend columns=3, 
	]
	\addplot table[x=k,y = t] {data/opt-T0Boltz-qr-fast-eim-S8-UNIa-trunc8-a-11-b110-restol1e-06-numTrain50-atrue0.5-btrue2-noise0.txt};
	\addlegendentry{ROM, bilevel, grad-free}
	\addplot table[x=k,y = t] {data/opt-T1Boltz-qr-fast-eim-S8-UNIa-trunc8-a-11-b110-restol1e-06-numTrain50-atrue0.5-btrue2-noise0.txt};
	\addlegendentry{ROM, bilevel, grad-based}
   	\addplot table[x=k,y = t] {data/optKKT-TBoltz-qr-fast-eim-S8-UNIa-trunc8-a-11-b110-restol1e-06-numTrain50-atrue0.5-btrue2-noise0.txt};
	\addlegendentry{ROM, single-level, grad-based}

	\nextgroupplot[
	ylabel={Reconstructed $a$},
	xlabel={},
	]

    \addplot table[x=k,y = a] {data/opt-T0Boltz-qr-fast-eim-S8-UNIa-trunc8-a-11-b110-restol1e-06-numTrain50-atrue0.5-btrue2-noise0.txt};
	\addplot table[x=k,y = a] {data/opt-T1Boltz-qr-fast-eim-S8-UNIa-trunc8-a-11-b110-restol1e-06-numTrain50-atrue0.5-btrue2-noise0.txt};
   	\addplot table[x=k,y = a] {data/optKKT-TBoltz-qr-fast-eim-S8-UNIa-trunc8-a-11-b110-restol1e-06-numTrain50-atrue0.5-btrue2-noise0.txt};
    \addplot[black, dashed, thick,domain=0:17,] {0.5};

	\nextgroupplot[
	ylabel={Reconstructed $b$},
	xlabel={},
	]
    \addplot table[x=k,y = b] {data/opt-T0Boltz-qr-fast-eim-S8-UNIa-trunc8-a-11-b110-restol1e-06-numTrain50-atrue0.5-btrue2-noise0.txt};
	\addplot table[x=k,y = b] {data/opt-T1Boltz-qr-fast-eim-S8-UNIa-trunc8-a-11-b110-restol1e-06-numTrain50-atrue0.5-btrue2-noise0.txt};
   	\addplot table[x=k,y = b] {data/optKKT-TBoltz-qr-fast-eim-S8-UNIa-trunc8-a-11-b110-restol1e-06-numTrain50-atrue0.5-btrue2-noise0.txt};
    \addplot[black, dashed, thick,domain=0:17,] {2};

	\nextgroupplot[
	ylabel={Data misfit},
	ymode=log,
    ymin = 3e-10, ymax = 3e-5,
	]
        \addplot[draw=none,name path=shadebottom,forget plot] coordinates {(10.5,3e-5) (16.5,3e-5)};
\addplot[draw=none,name path=shadetop,forget plot] coordinates {(10.5,3e-10) (16.5,3e-10)};
\addplot[fill=gray!20,draw=none,forget plot] fill between[of=shadebottom and shadetop];
    \addplot table[x=k,y = fval] {data/opt-T0Boltz-qr-fast-eim-S8-UNIa-trunc8-a-11-b110-restol1e-06-numTrain50-atrue0.5-btrue2-noise0.txt};
	\addplot table[x=k,y = fval] {data/opt-T1Boltz-qr-fast-eim-S8-UNIa-trunc8-a-11-b110-restol1e-06-numTrain50-atrue0.5-btrue2-noise0.txt};
   	\addplot table[x=k,y = fval] {data/optKKT-TBoltz-qr-fast-eim-S8-UNIa-trunc8-a-11-b110-restol1e-06-numTrain50-atrue0.5-btrue2-noise0.txt};

	\nextgroupplot[
	ylabel={Rel err $|a-a^*|/|a^*|$},
	ymode=log,
     ymin = 1e-4, ymax=3,
	]
            \addplot[draw=none,name path=shadebottom,forget plot] coordinates {(10.5,5e-6) (16.5,5e-6)};
\addplot[draw=none,name path=shadetop,forget plot] coordinates {(10.5,3) (16.5,3)};
\addplot[fill=gray!20,draw=none,forget plot] fill between[of=shadebottom and shadetop];
    \addplot table[x=k,y = relEa] {data/opt-T0Boltz-qr-fast-eim-S8-UNIa-trunc8-a-11-b110-restol1e-06-numTrain50-atrue0.5-btrue2-noise0.txt};
	\addplot table[x=k,y = relEa] {data/opt-T1Boltz-qr-fast-eim-S8-UNIa-trunc8-a-11-b110-restol1e-06-numTrain50-atrue0.5-btrue2-noise0.txt};
   	\addplot table[x=k,y = relEa] {data/optKKT-TBoltz-qr-fast-eim-S8-UNIa-trunc8-a-11-b110-restol1e-06-numTrain50-atrue0.5-btrue2-noise0.txt};

	\nextgroupplot[
	ylabel={Rel err $|b-b^*|/|b^*|$},
	ymode=log,
     ymin = 1e-4, ymax=3,
	]
\addplot[draw=none,name path=shadebottom,forget plot] coordinates {(10.5,9e-5) (16.5,9e-5)};
\addplot[draw=none,name path=shadetop,forget plot] coordinates {(10.5,3) (16.5,3)};
\addplot[fill=gray!20,draw=none,forget plot] fill between[of=shadebottom and shadetop];
    \addplot table[x=k,y = relEb] {data/opt-T0Boltz-qr-fast-eim-S8-UNIa-trunc8-a-11-b110-restol1e-06-numTrain50-atrue0.5-btrue2-noise0.txt};
	\addplot table[x=k,y = relEb] {data/opt-T1Boltz-qr-fast-eim-S8-UNIa-trunc8-a-11-b110-restol1e-06-numTrain50-atrue0.5-btrue2-noise0.txt};
   	\addplot table[x=k,y = relEb] {data/optKKT-TBoltz-qr-fast-eim-S8-UNIa-trunc8-a-11-b110-restol1e-06-numTrain50-atrue0.5-btrue2-noise0.txt};
	
	\end{groupplot}
	\node[black] at ($(group c2r1) + (0cm,2cm)$) {\pgfplotslegendfromname{grouplegend-ivp-Boltz-b110-positive-a-S8}};
	\end{tikzpicture}
	\caption{Inverse problem in \emph{transport-dominated} regime:  Comparison of ROM-based inverse problem solvers using 
    different optimization approaches.  The true parameters are $a^*=0.5$ and $b^*=2$, indicated by the black dashed lines. All solvers use initialization $\bmu = (-0.2,5)$. 
    The ROM uses the RB spaces (with different dimensions) trained in \Cref{fig:ROM-qr-fast-eim-S8-b110-data} in \Cref{sec:num-rom}.  The gray-shaded regions in subfigures refer to the RB dimensions for which a plateau is observed in the testing results (C) in \Cref{fig:ROM-qr-fast-eim-S8-b110-data}. 
	}\label{fig:ivp-Boltz-b110-positive-a-S8}
\end{figure}

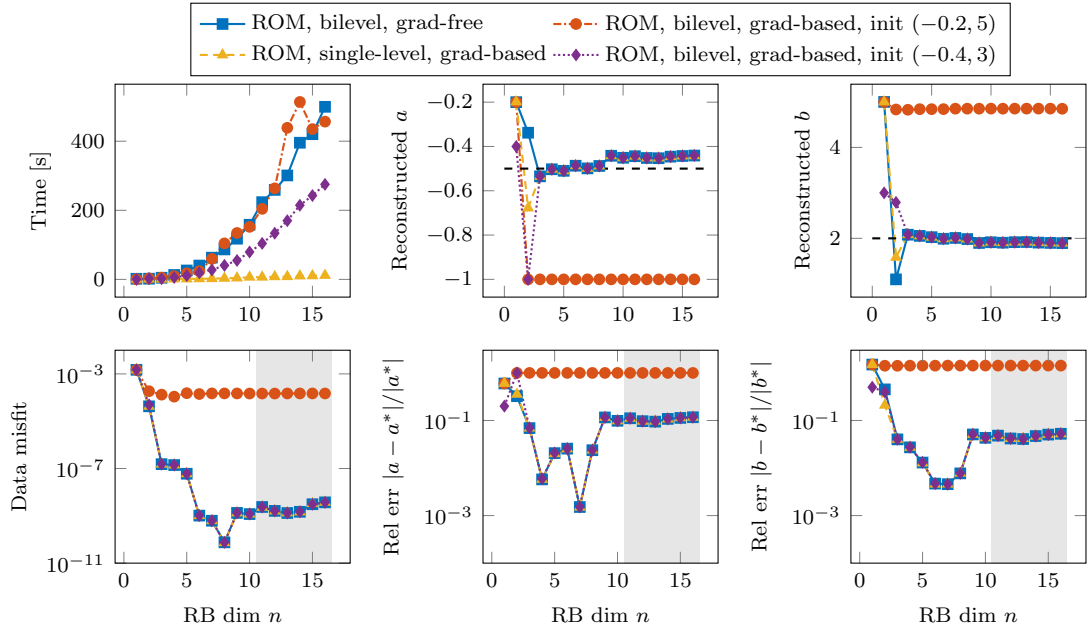
\begin{figure}[h]
	\centering
	\begin{tikzpicture}
	\begin{groupplot}[
	group style={
		group size=3 by 2,
		horizontal sep=50pt,
		vertical sep =20pt,
	},
	width=.32\textwidth,
	height=.3\textwidth,
	legend style={font=\scriptsize},
	tick label style={font=\scriptsize},
	label style={font=\scriptsize},
	xlabel = {RB dim $n$},
	legend style={font=\scriptsize},
	tick label style={font=\scriptsize},
	xmax = 18,
    cycle list name=siam,
    axis on top,
	]
	
	\nextgroupplot[
	xlabel={},
	ylabel={Time [s]},
	legend style = {nodes=right, legend to name=grouplegend-ivp-Boltz-b110-negative-a-S8},
	legend columns=2, 
	]
	\addplot table[x=k,y = t] {data/opt-T0Boltz-qr-fast-eim-S8-UNIa-trunc8-a-11-b110-restol1e-06-numTrain50-atrue-0.5-btrue2-noise0.txt};
	\addlegendentry{ROM, bilevel, grad-free}
	\addplot table[x=k,y = t] {data/opt-T1Boltz-qr-fast-eim-S8-UNIa-trunc8-a-11-b110-restol1e-06-numTrain50-atrue-0.5-btrue2-noise0-inita-0.2b5.txt};
	\addlegendentry{ROM, bilevel, grad-based, init $(-0.2,5)$}
   	\addplot table[x=k,y = t] {data/optKKT-TBoltz-qr-fast-eim-S8-UNIa-trunc8-a-11-b110-restol1e-06-numTrain50-atrue-0.5-btrue2-noise0.txt};
	\addlegendentry{ROM, single-level, grad-based}
	\addplot table[x=k,y = t] {data/opt-T1Boltz-qr-fast-eim-S8-UNIa-trunc8-a-11-b110-restol1e-06-numTrain50-atrue-0.5-btrue2-noise0-inita-0.4b3.txt};
	\addlegendentry{ROM, bilevel, grad-based, init $(-0.4,3)$}

	\nextgroupplot[
	ylabel={Reconstructed $a$},
	xlabel={},
	]

    \addplot table[x=k,y = a] {data/opt-T0Boltz-qr-fast-eim-S8-UNIa-trunc8-a-11-b110-restol1e-06-numTrain50-atrue-0.5-btrue2-noise0.txt};
	\addplot table[x=k,y = a] {data/opt-T1Boltz-qr-fast-eim-S8-UNIa-trunc8-a-11-b110-restol1e-06-numTrain50-atrue-0.5-btrue2-noise0-inita-0.2b5.txt};
   	\addplot table[x=k,y = a] {data/optKKT-TBoltz-qr-fast-eim-S8-UNIa-trunc8-a-11-b110-restol1e-06-numTrain50-atrue-0.5-btrue2-noise0.txt};
	\addplot table[x=k,y = a] {data/opt-T1Boltz-qr-fast-eim-S8-UNIa-trunc8-a-11-b110-restol1e-06-numTrain50-atrue-0.5-btrue2-noise0-inita-0.4b3.txt};
    \addplot[black, dashed, thick,domain=0:17,] {-0.5};

	\nextgroupplot[
	ylabel={Reconstructed $b$},
	xlabel={},
	]
    \addplot table[x=k,y = b] {data/opt-T0Boltz-qr-fast-eim-S8-UNIa-trunc8-a-11-b110-restol1e-06-numTrain50-atrue-0.5-btrue2-noise0.txt};
	\addplot table[x=k,y = b] {data/opt-T1Boltz-qr-fast-eim-S8-UNIa-trunc8-a-11-b110-restol1e-06-numTrain50-atrue-0.5-btrue2-noise0-inita-0.2b5.txt};
   	\addplot table[x=k,y = b] {data/optKKT-TBoltz-qr-fast-eim-S8-UNIa-trunc8-a-11-b110-restol1e-06-numTrain50-atrue-0.5-btrue2-noise0.txt};
	\addplot table[x=k,y = b] {data/opt-T1Boltz-qr-fast-eim-S8-UNIa-trunc8-a-11-b110-restol1e-06-numTrain50-atrue-0.5-btrue2-noise0-inita-0.4b3.txt};
    \addplot[black, dashed, thick,domain=0:17,] {2};

	\nextgroupplot[
	ylabel={Data misfit},
	ymode=log,
    ymin = 1e-11, ymax = 1e-2,
	]
 \addplot[draw=none,name path=shadebottom,forget plot] coordinates {(10.5,1e-11) (16.5,1e-11)};
\addplot[draw=none,name path=shadetop,forget plot] coordinates {(10.5,3) (16.5,3)};
\addplot[fill=gray!20,draw=none,forget plot] fill between[of=shadebottom and shadetop];
    \addplot table[x=k,y = fval] {data/opt-T0Boltz-qr-fast-eim-S8-UNIa-trunc8-a-11-b110-restol1e-06-numTrain50-atrue-0.5-btrue2-noise0.txt};
	\addplot table[x=k,y = fval] {data/opt-T1Boltz-qr-fast-eim-S8-UNIa-trunc8-a-11-b110-restol1e-06-numTrain50-atrue-0.5-btrue2-noise0-inita-0.2b5.txt};
   	\addplot table[x=k,y = fval] {data/optKKT-TBoltz-qr-fast-eim-S8-UNIa-trunc8-a-11-b110-restol1e-06-numTrain50-atrue-0.5-btrue2-noise0.txt};
	\addplot table[x=k,y = fval] {data/opt-T1Boltz-qr-fast-eim-S8-UNIa-trunc8-a-11-b110-restol1e-06-numTrain50-atrue-0.5-btrue2-noise0-inita-0.4b3.txt};

	\nextgroupplot[
	ylabel={Rel err $|a-a^*|/|a^*|$},
	ymode=log,
      ymin = 1e-4, ymax=3,
	]
         \addplot[draw=none,name path=shadebottom,forget plot] coordinates {(10.5,1e-11) (16.5,1e-11)};
\addplot[draw=none,name path=shadetop,forget plot] coordinates {(10.5,3) (16.5,3)};
\addplot[fill=gray!20,draw=none,forget plot] fill between[of=shadebottom and shadetop];
    \addplot table[x=k,y = relEa] {data/opt-T0Boltz-qr-fast-eim-S8-UNIa-trunc8-a-11-b110-restol1e-06-numTrain50-atrue-0.5-btrue2-noise0.txt};
	\addplot table[x=k,y = relEa] {data/opt-T1Boltz-qr-fast-eim-S8-UNIa-trunc8-a-11-b110-restol1e-06-numTrain50-atrue-0.5-btrue2-noise0-inita-0.2b5.txt};
   	\addplot table[x=k,y = relEa] {data/optKKT-TBoltz-qr-fast-eim-S8-UNIa-trunc8-a-11-b110-restol1e-06-numTrain50-atrue-0.5-btrue2-noise0.txt};
	\addplot table[x=k,y = relEa] {data/opt-T1Boltz-qr-fast-eim-S8-UNIa-trunc8-a-11-b110-restol1e-06-numTrain50-atrue-0.5-btrue2-noise0-inita-0.4b3.txt};

	\nextgroupplot[
	ylabel={Rel err $|b-b^*|/|b^*|$},
	ymode=log,
     ymin = 1e-4, ymax=3,
	]
     \addplot[draw=none,name path=shadebottom,forget plot] coordinates {(10.5,1e-11) (16.5,1e-11)};
\addplot[draw=none,name path=shadetop,forget plot] coordinates {(10.5,3) (16.5,3)};
\addplot[fill=gray!20,draw=none,forget plot] fill between[of=shadebottom and shadetop];
    \addplot table[x=k,y = relEb] {data/opt-T0Boltz-qr-fast-eim-S8-UNIa-trunc8-a-11-b110-restol1e-06-numTrain50-atrue-0.5-btrue2-noise0.txt};
	\addplot table[x=k,y = relEb] {data/opt-T1Boltz-qr-fast-eim-S8-UNIa-trunc8-a-11-b110-restol1e-06-numTrain50-atrue-0.5-btrue2-noise0-inita-0.2b5.txt};
   	\addplot table[x=k,y = relEb] {data/optKKT-TBoltz-qr-fast-eim-S8-UNIa-trunc8-a-11-b110-restol1e-06-numTrain50-atrue-0.5-btrue2-noise0.txt};
	\addplot table[x=k,y = relEb] {data/opt-T1Boltz-qr-fast-eim-S8-UNIa-trunc8-a-11-b110-restol1e-06-numTrain50-atrue-0.5-btrue2-noise0-inita-0.4b3.txt};
	
	\end{groupplot}
	\node[black] at ($(group c2r1) + (0cm,2cm)$) {\pgfplotslegendfromname{grouplegend-ivp-Boltz-b110-negative-a-S8}};
	\end{tikzpicture}
	\caption{Inverse problem in \emph{transport-dominated} regime: Comparison of ROM-based inverse problem solvers  using 
    different optimization approaches. The true parameters are $a^*=-0.5$ and $b^*=2$, indicated by the black dashed lines. The gradient-free bilevel and single-level KKT solvers use initialization $\bmu = (-0.2,5)$. For the gradient-based bilevel solver, we report two initializations, $(-0.2,5)$ and $(-0.4,3)$. The ROM uses the RB spaces (with different dimensions) trained in \Cref{fig:ROM-qr-fast-eim-S8-b110-data} in \Cref{sec:num-rom}. The gray-shaded regions in subfigures refer to the RB dimensions for which a plateau is observed in the testing results (C) in \Cref{fig:ROM-qr-fast-eim-S8-b110-data}. 
	}\label{fig:ivp-Boltz-b110-negative-a-S8}
\end{figure}

\medskip
\noindent {\bf Inverse problem in transport-dominated regime.} By using the ROM trained for the relatively transport-dominated regime $\mM=[-1,1]\times[1,10]$,
we now test the three optimization approaches to solve the inverse problem, when the true parameter $\bmu^*=(a^*, b^*)$ is one of the following: $(0.5, 2)$ and $(-0.5, 2)$. %
For all cases, the initial guess of the parameter is 
set with $(-0.2, 5)$ unless otherwise stated, and the results are  shown  in \Cref{fig:ivp-Boltz-b110-positive-a-S8,fig:ivp-Boltz-b110-negative-a-S8}.  As before, the gray-shaded region corresponds to the RB dimensions (i.e., $n\geq 11$) for which the plateau effect is observed in the testing errors of our ROM as in \Cref{fig:ROM-qr-fast-eim-S8-b110-data}.

For the case with $a^*=0.5$ and $b^*=2$, the data misfit and the relative reconstruction errors in \Cref{fig:ivp-Boltz-b110-positive-a-S8} show an overall decreasing trend as the ROM solvers improve their resolution until around $n=9$ (that is, slightly before the onset of the plateau), after which they may increase, just as in the collision-dominated regime.  The best reconstruction is obtained at $n=9$, with relative reconstruction errors of approximately $0.1\%$. In addition, the single-level optimization method achieves a relative error of approximately $0.01\%$ for $b$. The improved accuracy in reconstruction  compared with that in the collision-dominated regime is related to the greater parameter sensitivity of the solutions in the transport-dominated regime. Although this stronger sensitivity requires slightly larger RB spaces during ROM training to capture the solution manifold accurately, it makes the inverse problem easier because the observations and solutions respond more strongly to changes in the parameters.

For the case with $a^*=-0.5$ and $b^*=2$, similar observations as above can be made from  \Cref{fig:ivp-Boltz-b110-negative-a-S8}, with the best reconstruction obtained at $n=8$ and the associated relative reconstruction errors being  approximately $0.1\%$. The only exception is the gradient-based bilevel approach, which becomes trapped at the boundary of the parameter domain. 
The same phenomenon is observed when the FOM is used instead of the ROM to solve the inverse problem. The optimization requires $18361$ seconds to converge, but the reconstructed parameter becomes trapped on the boundary at $(-1, 4.85)$ with a data misfit of approximately $1.5\times 10^{-4}$. The performance of the gradient-based bilevel ROM approach improves when a closer initial guess $\bmu=(-0.4,3)$ is used, as shown in the purple dotted lines in \Cref{fig:ivp-Boltz-b110-negative-a-S8}.

\medskip
\noindent{\bf Comparison with FOM-based inverse solvers.}
Finally, we want to compare the three ROM-based inverse solvers  with that based on FOM regarding  both the computational efficiency and the reconstruction accuracy.  In the transport-dominated regime, the FOM requires less computational time than that in the collision-dominated regime, making it feasible to solve the inverse problem directly using the FOM. The FOM-based inverse solver is analogous to the bilevel formulation \eqref{eq:IP-Boltz-ROM-bilevel}, except that the inner constrained optimization in \Cref{alg:bilevel} is replaced by the iterative FOM solver described in \Cref{sec:FOM}. For the outer optimization, we use the quasi-Newton algorithm in \texttt{fminunc} with the same settings as those used for the gradient-based bilevel ROM-based optimization approach. We focus on the case with the true parameter $\bmu^*=(0.5, 2)$.
The FOM-based inverse solver requires $12085$ seconds in total, with $105$ FOM solves, and $29$ optimization iterations. The relative errors in the reconstructed parameters $a$ and $b$ are $1.4\times 10^{-5}$ and $8.1\times 10^{-6}$, respectively, and the final data misfit is $2\times 10^{-16}$. This yields approximately two additional orders of magnitude in reconstruction accuracy compared with the ROM-based approaches, but at a substantially higher computational cost. With the RB dimension $n=9$, the FOM-based method is approximately $90$ times slower than the gradient-free bilevel ROM approach, $190$ times slower than the gradient-based bilevel ROM approach, and $3000$ times slower than the gradient-based single-level approach.

We further consider the FOM-based inverse solver with  an early termination when the data misfit falls below $10^{-9}$, which is comparable to the misfit associated with the best reconstruction shown in \Cref{fig:ivp-Boltz-b110-positive-a-S8}. In this case, the FOM-based method requires $11931$ seconds, $93$ FOM solves,  and $25$ optimization iterations. The relative reconstruction errors in $a$ and $b$ are $1.1\times 10^{-2}$ and $6.4\times 10^{-3}$, respectively, and the final data misfit is $1.4\times10^{-10}$. These reconstruction errors are larger than those obtained by the three ROM-based approaches, 
while the computational time remains hundreds to thousands of times greater.
The improved reconstruction obtained with the ROM, compared with the FOM as the forward solver, may result from restricting the inverse problem to a lower-dimensional reduced space, leading to  a better optimization landscape and partially mitigating the ill-posedness of the problem. 
A similar phenomenon was also observed in \cite{borcea2024data}.

\section{Conclusion and discussion}
\label{sec:con}

In this paper, we present a model order reduction framework that provides a reduced-order model (ROM) for solving the steady-state parametric Boltzmann equation. Significant computational efficiency is achieved via low-dimensional representation/approximation spaces built greedily offline for the parameter-induced solution manifold, and precomputation associated with the nonlocal and nonlinear collision operator based on its quadratic structure and a separable EIM kernel approximation. Physical structure  such as mass conservation  is also explicitly enforced.  
 The ROM and the associated reduced-basis space are further applied to a thermally-driven inverse problem  to infer the parameters in the collision kernel from the observed macroscopic temperature data.  Here, both a direct bilevel optimization formulation and a single-level nonlinear constrained optimization formulation are considered to solve the inverse problem.

Though the methodologies are applicable to the general Boltzmann equation,  they are described and numerically demonstrated for the case with one-dimensional physical space and two-dimensional velocity space.  The resulting ROM achieves more than three orders of magnitude acceleration as a surrogate solver  
while maintaining an accuracy of approximately $10^{-4}$, compared to the full-order model.  
Although the inverse problem is highly ill-posed and computationally expensive if standard forward solvers are used,  the proposed ROM and the reduced representation of the solution  make the original optimization approaches to solve the inverse problem computationally tractable, reducing the solve time to less than $10$ seconds in representative cases. In addition, some numerical evidence shows  that working in the reduced space may improve the reconstruction quality by possibly  providing a more favorable optimization landscape. %

Since the Boltzmann equation and its inverse problems have attracted increasing attention because of their broad applications and mathematical significance, %
the present work can be viewed as a first step  
toward developing ROM-based approaches to reduce the computational complexity of the corresponding forward and inverse problems.
Natural extensions include higher-dimensional phase-space settings, more general collision kernels (e.g., with angular dependence), and richer experimental configurations and observations---such as combining Fourier and shear-driven Couette flows---to improve parameter identifiability.
Another issue is the reconstruction from noisy data. Given that the inverse problem is very ill-posed, numerical experiments are not presented when our proposed algorithms are applied to noisy data. 
Extending the proposed framework to achieve stable reconstruction from noisy data remains an important direction for future work.

\section*{Acknowledgment}
All authors acknowledge the support of the ICERM workshop ``Empowering a Diverse Computational Mathematics Research Community" on July 22 - August 2, 2024, where this work was initiated.

ST is partially supported by National Science Foundation under Grant No.~DMS-2529292. JH is partially supported by National Science Foundation under Grant No.~DMS-2409858. FL is partially supported by Air Force Office of Scientific Research under Grant No.~FA9550-26-1-0003. ZS is partially supported by SIAM Postdoctoral Support Program. YY is partially supported by National Science Foundation under Grant No.~DMS-2409855 and DMS-2540324, and Office of Naval Research under Award No.~N00014-24-1-2088.

\section*{Declaration of generative AI and AI-assisted technologies in the manuscript preparation process}
During the preparation of this work, some of the authors used AI-based language tools (e.g., ChatGPT) in order to improve phrasing. After using this tool, the authors reviewed and edited the content as needed and take full responsibility for the content of the published article.

\bibliography{ROM,hu_bibtex}
	\bibliographystyle{plain}

\end{document}